\documentclass[11pt,a4paper, twoside]{article}
\usepackage[notref,notcite]{showkeys}
\usepackage{amsmath,amssymb,amsthm,amsfonts,mathrsfs,amscd,environ}
\usepackage{dsfont}
\usepackage{latexsym,enumerate,color,geometry}
\usepackage{xcolor}
\usepackage{verbatim,fancyhdr}
 \NewEnviron{ews}{%
\begin{equation}\begin{split}
  \BODY
\end{split}\end{equation}
}

\NewEnviron{ews*}{%
\begin{equation*}\begin{split}
  \BODY
\end{split}\end{equation*}
}

\def\beg{\begin}
\def\bequ{\begin{equation}}
\def\enqu{\end{equation}}
\def\bes{\begin{split}}
\def\ens{\end{split}}
\def\bews{\begin{ews}}
\def\beqn{\begin{eqnarray}}
\def\enqn{\end{eqnarray}}
\def\beq*{\begin{equation*}}
\def\enq*{\end{equation*}}
\def\bqn*{\begin{eqnarray*}}
\def\eqn*{\end{eqnarray*}}
\def\bary{\begin{array}}
\def\eary{\end{array}}
\def\bpma{\begin{pmatrix}}
\def\epma{\end{pmatrix}}
\def\bvma{\begin{Vmatrix}}
\def\evma{\end{Vmatrix}}

 \numberwithin{equation}{section}

\def\al{\alpha}
\def\be{\beta}
\def\ga{\gamma}
\def\de{\delta}
\def\ep{\epsilon}

\def\ze{\zeta}
\def\et{\eta}
\def\th{\theta}

\def\la{\lambda}
\def\rh{\rho}

\def\si{\sigma}
\def\ta{\tau}

\def\ps{\psi}
\def\om{\omega}

\def\La{\Lambda}
\def\Ga{\Gamma}

\def\De{\Delta}

\def\Ps{\Psi}
\def\Om{\Omega}

\def\R{\mathbb R}
\def\P{\mathbb P}
\def\E{\mathbb E}

\def\N{\mathbb N}

\def\SS{\mathbb S}

\def\sE{\mathscr E}
\def\sF{\mathscr F}
\def\sD{\mathscr D}

\def\sB{\mathscr B}

\def\sG{\mathscr G}

\def\cO{\mathcal O}

\def\d{\mathrm{d}}

\def\ff{\frac}
\def\ra{\rightarrow}
\def\nn{\nabla}
\def\pp{\partial}
\def\<{\langle}
\def\>{\rangle}
\def\sq{\sqrt}
\def\tld{\tilde}
\def\we{\wedge}
\def\1{\mathds{1}}

\def\trac{\mathrm{tr}}

\def\div{\mathrm{div}}

\def\supp{\displaystyle{\mathrm{supp}}}

\allowdisplaybreaks

\title{{\bf On  invariant probability measures of regime-switching diffusion processes with singular drifts}
}

\author{
{\bf Shao-Qin Zhang }\\
\footnotesize{School of Statistics and Mathematics, Central University of Finance and Economics, Beijing 100081, China}\\
\footnotesize{Email: zhangsq@cufe.edu.cn}\\
}

\begin{document}
\maketitle

\begin{abstract}
We introduce integrability conditions involving a nice reference probability measure and the $Q$-matrix to  study the existence  of the invariant probability measures of regime-switching diffusion processes. Regularities and the uniqueness  of the invariant probability density w.r.t the nice reference probability measure are also considered. Moreover, we study the $L^1$-uniqueness of the semigroup generated by regime-switching diffusion processes. Since the generator of the regime-switching diffusion process is a weakly coupled elliptic system,  the $L^1$-uniqueness of the extension of  weakly coupled elliptic systems is obtained.
\end{abstract}\noindent

\noindent


\vskip 2cm

\section{Introduction}
Let  $\SS=\{1,2,3,\cdots,N\}$, and let $(\Om,\sF,\P)$ be a completed probability space. A regime-switching diffusion process is a two components process $\{(X_t,\Lambda_t)\}_{t\geq 0}$ described by 
\bequ\label{equ_b1}
\d X_t=b_{\Lambda_t}(X_t)\d t+\sq 2 \si_{\Lambda_t}(X_t)\d W_t,
\enqu
and
\bequ\label{equ_b2}
\P\left(\Lambda_{t+\De t}=j~|~\Lambda_t=i,~(X_s,\Lambda_s),s\leq t\right)
=\left\{\beg{array}{lr}
q_{ij}(X_t)\De t+o(\De t),& i\neq j,\\
1+q_{ii}(X_t)\De t+o(\De t),& i=j,
\end{array}
\right.
\enqu
where $\{W_t\}_{t\geq 0}$ is a Brownian motion on $\R^d$ w.r.t a right continuous completed reference $\{\sF_t\}_{t\geq 0}$ and 
$$b:\R^d\times \SS \ra \R^d,~~\si:\R^d\times \SS \ra \R^d\otimes\R^d,~~q_{ij}:\R^d\ra \R,$$
are  measurable functions and 
$$q_{ij}(x)\geq 0,~i\neq j,~~~~~\sum_{j\neq i}q_{ij}(x)\leq -q_{ii}(x).$$
The matrix $Q(x)=(q_{ij}(x))_{1\leq i,j\leq N}$ is called $Q$-matrix, see \cite{ChenB}. If $Q(x)$ is independent of $x$ and $\Lambda_t$ is independent of $\{W_t\}_{t\geq 0}$, $\{(X_t,\Lambda_t)\}_{t\geq 0}$ is called a state-independent RSDP, otherwise called a state-dependence RSDP. Let 
$$q_i(x)=\sum_{j\neq i}q_{ij}(x),~i\in\SS.$$ If $Q(x)$ is conservative, i.e. $-q_{ii}(x)= q_i(x)$, $x\in\R^d$, then as in \cite[Chapter II-2.1]{Skoro} or \cite{GAM,YinZ09,Shao15}, we can represent $\{(X_t,\Lambda_t)\}_{t\geq 0}$ in the form of a system of stochastic differential equations(SDEs for short) driven by $\{W_t\}$ and a Poisson random measure. Precisely, for each $x\in\R^d$, $\{\Ga_{ij}(x):~i,~j\in\SS\}$ is a family of disjoint  intervals on $[0,\infty)$ constructed as follows
\beg{align*}
&\Ga_{12}(x)=[0,q_{12}(x)),~\Ga_{13}(x)=[q_{12}(x),q_{12}(x)+q_{13}(x)),\cdots\\
&\Ga_{21}(x)=[q_1(x),q_1(x)+q_{21}(x)),~\Ga_{23}(x)=[q_1(x)+q_{21}(x),q_1(x)+q_{21}(x)+q_{23}(x)),\cdots\\
&\Ga_{31}(x)=[q_1(x)+q_2(x),q_1(x)+q_{2}(x)+q_{13}(x)),\cdots\\
&\qquad\cdots
\end{align*}
We set $\Ga_{ii}(x)=\emptyset$ and $\Ga_{ij}(x)=\emptyset$ if $q_{ij}(x)=0$ for $i\neq j$. Define a function $h:\R^d\times\SS\times [0,\infty)\ra\R$  as follows
\beq*
h(x,i,z)=\sum_{j\in\SS}(j-i)\1_{\Ga_{ij}(x)}(z)=\left\{\beg{array}{cl}
j-i,& \mbox{if~} z \in \Ga_{ij}(x),\\
0,& \mbox{otherwise}.
\end{array}
\right.
\enq*
Let $N(\d z,\d t)$ be a Poisson random measure with intensity $\d z\d t$ and independent of the Brownian motion $\{W_t\}_{t\geq 0}$. Then we turn to consider the following equation
\bequ\label{main_equ}
\beg{cases}
\d X_t = b_{\Lambda_t}(X_t)\d t+\sq 2\si_{\Lambda_t}(X_t)\d W_t,\\
\d \Lambda_t = \int_{0}^\infty h(X_{t-},\Lambda_{t-},z)N(\d z,\d t).
\end{cases}
\enqu
The first component $X_t$ can be viewed as a hybrid process from the diffusion $X^i_t$ in each state in $\SS$:
\bequ\label{Xi}
\d X_t^i=b_i(X_t^i)\d t+\sq 2\si_i(X_t^i)\d W_t,~i\in\SS.
\enqu
In \cite{ZhangS}, we have study the existence, uniqueness  and strong Feller property of \eqref{main_equ}, where drifts $b_i$, $i=1,\cdots, N$ can be singular. In this paper, we shall study the existence and uniqueness of invariant probability measures of \eqref{main_equ} with non-regular drifts.

Recently, in term of a nice reference probability measure, \cite{Wang17,Wang18} introduce integrability conditions to ensure the existence of invariant probability measures for SDEs with singular or path dependent drifts. We shall extract similar results in \cite{Wang17,Wang18} for the semigroup corresponding to \eqref{main_equ}. However, for this type of semigroups generated by \eqref{main_equ}, two obstacles have to be concern about. One is that, as pointed out by \cite{Shao}, to have an invariant probability measure, it is not necessary for each diffusion $X_t^i$ to process an invariant probability measure. The other one is the heavy tail phenomenon found in \cite{DYao}.  Hence the integrability of invariant measures of \eqref{main_equ} may be worse than the invariant probability measures of some $X_t^i$. So new integrability conditions (see \eqref{ewZ}, \eqref{inequ_MM} and Remark \ref{re_exa}) are needed to study to invariant probability measure of \eqref{main_equ}.  

Let $a_k(x)=\si_k(x)\si^*_k(x)=(a_k^{i,j}(x))_{1\leq i,j\leq d}$  and $\div(a_j(x))_l=\sum_{k=1}^d\pp_k a_j^{kl}(x)$.  Let $V\in C^2(\R^d)$ such that $\mu_V(\d x)=e^{V(x)}\d x$ is a probability measure.  In this paper, we consider the drifts that process the following form
$$b_k(x)\equiv Z_k^0(x)+\sq 2 \si_k(x)Z(x)\equiv a_k(x)\nn V(x)+\div(a_k(x))+\sq 2\si_k(x)Z_k(x).$$
Then $\si_kZ_k$ is the singular part of the drift $b_k$. For the diffusion part,  we define  differential  operators $L^0_k$ and $L_k^Z$, $k=1,\cdots,N$   as follows
\beg{align*}
(L^0_kf)(x)&=e^{-V(x)}\div\left(e^{V(x)}a_k(x)\nn f(x)\right)\nonumber\\
&=\trac\left(a_k(x)\nn^2 f(x)\right)+\left\<a_k(x)\nn V(x),\nn f(x)\right\>\\
&\qquad +\Big\<\div(a_k(x)),\nn f(x)\Big\>\nonumber\\
&\equiv \trac\left(a_k(x)\nn^2 f(x)\right)+\Big\<Z_k^0(x),\nn f(x)\Big\>,\nonumber\\
(L_k^Z f)(x)&=(L^0_k f)(x)+\sq 2\<\si_k(x)Z_k(x),\nn f(x)\>,~~f\in C_c(\R^d).
\end{align*}
Let $P_t$ be the semigroup associated with \eqref{main_equ}.   Then  the generator of $P_t$, denoted by $L$, is of the following form
\beg{align*}
(L f)(x,k)&= L_k^Zf_k(x)+\sum_{j=1}^Nq_{kj}(x)f_j(x),~f\in C_c^2(\R^d\times\SS).
\end{align*}
The  operator $L$ sometimes is called weakly coupled elliptic system, see \cite{ChenZ}.

Let $\mu_{\pi,V}$ be a probability measure on $\R^d\times\SS$:
$$\mu_{\pi,V}(f)=\sum_{k=1}^N\int_{\R^d}f_k(x) \pi_ke^{V(x)}\d x.$$
Denote by $H_{\si_j}^{1,2}(\R^d)$ the Sobolev space which is the completion of $C^2_c(\R^d)$ under the norm
$$\left(\int_{\R^d}\left(|\si_j^*(x)\nn u(x)|^2+|u(x)|^2\right)\mu_V(\d x)\right)^{\ff 1 2},~u\in C^2_c(\R^d)$$ 
and by $W_{\si_j}^{1,2}(\R^d)$  the Sobolev space defined as follows
$$W_{\si_j}^{1,2}(\R^d)=\left\{u\in W^{1,2}_{loc}(\R^d)~\Big|~\int_{\R^d}\left(|\si_j^*(x)\nn u(x)|^2+|u(x)|^2\right)\mu_V(\d x)<\infty\right\}.$$ 
We assume that the coefficients of $L$ satisfy the following assumptions

\beg{enumerate}[(H1)]
\item $V\in C^2(\R^d)$ such that $e^{V}\d x$ is a probability measure, and $\si_k^{ij}\in C^2(\R^d)$. For all $k\in\SS$, 
 \beg{align}\label{H=W}
H_{\si_k}^{1,2}(\R^d)=W_{\si_k}^{1,2}(\R^d),
\end{align}
and there is $\la_k>0$ such that \label{H1}
\beg{align}
a_k(x)\geq \la_k,~x\in\R^d.
\end{align}

\item For $x\in\R^d$, $Q(x)=(q_{kj}(x))_{1\leq k,j\leq N}$ is a conservative $Q$-matrix with
$$C_Q:=\sup_{x\in\R^d,~k\in\SS} q_{k}(x)<\infty.$$
The matrix $Q$  is fully coupled on $\R^d$, i.e. the set $\SS$ can not be split into two disjoint non-empty $\SS_1$ and $\SS_2$ such that 
$$q_{l_1l_2}(\cdot)=0,\quad \mbox{a.e.},\qquad l_1\in\SS_1, l_2\in\SS_2.$$
There is an $x\in\R^d$ such that $Q(x)$ being an irreducible $Q$-Matrix with an invariant probability measure $\pi=(\pi_1,\cdots,\pi_N)$.  \label{H2}


\item  There are $\ga_1,\cdots,\ga_N\in (0,\infty)$ and $\be_1,\dots,\be_N\in [0,\infty)$ such that \label{H3}
$$\mu_V(f^2\log f^2)\leq \ga_k\mu_V(|\si_k^*\nn f|^2)+\be_k,~f\in H_{\si_k}^{1,2}(\R^d),~\mu_V(f^2)=1~k\in\SS.$$
\end{enumerate}

\beg{rem}\label{rem_logSob}
The condition (H\ref{H1}) implies that for all $i\in\SS$, $L^0_i$ generates a unique non-explosive Markov semigroup.

Let
$$\sE^0(f,g)=\sum_{k=1}^N\pi_k\int_{\R^d}\<\si_k^*\nn f_k(x),\si_k^*\nn g_k(x)\>\mu_{V}(\d x),~f_k,g_k\in H^{1,2}_{\si_k}(\R^d).$$
Then $\sE^0$ is a Dirichlet form with $\sD(\sE^0)=\{f=(f_1,\cdots,f_N)~|~f_i\in H_{\si_i}^{1,2}(\R^d),~i\in\SS\}$. The condition (H\ref{H3}) yields that the following defective log-Sobolev inequality holds. For $f_1\in H_{\si_1}^{1,2}(\R^d),~\cdots,~f_N\in H_{\si_N}^{1,2}(\R^d)$, with $\displaystyle\sum_{k=1}^N \pi_k\mu_V(f_k^2)=1$,
\beg{align}\label{log_N}
\mu_{\pi,V}(f^2\log f^2 )&=\sum_{k=1}^N\pi_k\mu_{V}(f_k^2\log f_k^2)\\
&\leq (\max_{k}\ga_k)\sE^{0}(f,f)+\sum_{k=1}^N\be_k\pi_k \mu_V(f_k^2)\nonumber\\
&\quad+\sum_{k=1}^N\pi_k\mu_V(f_k^2)\log\mu_V(f_k^2)\nonumber\\
&\leq \ga\sE^0(f,f)+\be,
\end{align}
with 
\end{rem}

\beg{rem}
The matrix $Q$ is fully coupled in the sense of \cite{Sw}. It was proved in \cite[Proposition 4.1.]{ChenZ96} that $Q$ is fully coupled  on a domain $D$ of $\R^d$ is equivalent to $Q$ is irreducible on $D$: for any distinct $k,l\in\SS$, there exist $k_0,k_1,\cdots,k_r$  in 􏰐$\SS$ with $k_i\neq k_{i+1}$, $k_0=􏰉k$ and $k_r=􏰉l$ such that $\{􏰐x\in D|~ q_{k,k+1}(x)\neq 0\}$  has positive Lebesgue measure for $i=􏰉0, 1,\cdots, r-1$.  
\end{rem}

\beg{thm}\label{thm_eu}
Assume (H\ref{H1})-(H\ref{H3}), and 
\bequ\label{infq}
\inf_{x\in \R^d}q_i(x)>0,~i\in\SS.
\enqu
Let $\bar Q=(\bar q_{ij})_{1\leq i,j\leq N}$ be a matrix with $\bar q_{ij}=\sup_{x\in\R^d}q_{ij}(x)$. We suppose in addition that  there exist positive vector $v=(v_1,v_2,\cdots,v_N)$ and positive constants $w_1,w_2,\cdots, w_N$ such that 
\beg{align}
\sup_{1\leq k\leq N}\mu_V(e^{w_k|Z_k|^2})&<\infty,\label{ewZ}\\
-(K+\bar Q)v &\geq \1, \label{inequ_MM}
\end{align}
where 
$$K=\mathrm{Diag}\left(\ff 1 {2 w_1}-\ff 2 {\ga_1},\ff 1 {2 w_2}-\ff 2 {\ga_2},\cdots,\ff 1 {2 w_N}-\ff 2 {\ga_N}\right).$$
Then $P_t$ has a unique invariant probability measure $\mu$ such that $\mu\ll \mu_{\pi,V}$. Moreover, the density  $\rh=\ff {\d \mu} {\d \mu_{\pi,V}}$ is positive on $\R^d\times \SS$,   and for all $p>0$ and $k\in\SS$, $\rh_k(\cdot)\equiv \rh(\cdot,k)\in W^{1,p}_{loc}(\R^d)$,  and $\sq\rh_k \in H_{\si_k}^{1,2}(\R^d)$. 
\end{thm}

\beg{rem}
A square matrix $A$  is called an M-Matrix if $A = sI -B$ with some $B \geq 0$ and $s\geq r(B)$, where by $B\geq 0$ we mean that all elements of $B$ are non-negative, $I$ is the  identity matrix, and $r(B)$ is the spectral radius of $B$. When $s>r(B)$, $A$ is called a nonsingular M-matrix. The inequality \eqref{inequ_MM} is equivalent to that $-(K+\bar Q)$ is a nonsingular M-Matrix, see \cite[Proposition 2.4]{Shao} or \cite{BP} for instance.

In our theorem, we indeed  prove   the uniqueness  in the  sense that the set 
$$\left\{\rh\mu_{\pi,V}~|~\rh\geq 0,~\mu_{\pi,V}(\rh)=1,~\mu_{\pi,V}(\rh P_t f)=\mu_{\pi,V}(\rh f)~\mbox{for}~t\geq 0,~f\in\sB(\R^d\times\SS)\right\}$$ 
contains only one element. 

According to  \cite[Theorem 2.1]{Wang17}, (H\ref{H1}) and   \eqref{ewZ} yields that for each $i$, the SDE \eqref{Xi} with generator  $L^Z_i$ has a non-explosive strong solution. Combining with that $\SS$ is finite and  $\cup_{x\in\R^d,i\in\SS}\supp(h(x,i,\cdot))$ is a bounded set of $[0,\infty)$, it is easy to see, for instance \cite{CH,ZhangS}, that \eqref{main_equ} has a non-explosive strong solution. 
\end{rem}

\beg{rem}\label{re_exa}
In conditions \eqref{ewZ} and \eqref{inequ_MM}, some $w_i$ can be small such that $w_i<\ff {\ga_i} 4 $, so not all  diffusion processes $X_t^i$ satisfy the integrability conditions in \cite[Thoerem 2.2]{Wang17} and \cite[Theorem 1.1(2) or Theorem 5.2]{Wang18}.  A concrete example is presented as follows.
\end{rem}
To illustrate that Theorem \ref{thm_eu} can be applied to the regime-switching diffusion process that not all the diffusions in all the environments has an invariant probability measure, we present the following example
\beg{exa}
Let $N=2$, $d=1$, $b_1(x)=-x$, $b_2(x)=-x+\sq 2\de x$, $\si_1(x)=\si_2(x)=1$, and 
$$Q(x)=\left(\begin{array}{cc} -a  & a\\ b  &-b \end{array} \right)+\left(\begin{array}{cc} -a(x)  & a(x)\\ b(x)  &-b(x) \end{array} \right)$$
with $a>0$, $b>0$, and $|a(x)|\leq \th a$,~$|b(x)|\leq \th b$ for some $\th\in [0,1)$. Then $V(x)=c-\ff 1 2 x^2$, $\ga_1=\ga_2=2$, $\be_1=\be_2=0$, $Z_1=0$ and  $Z_2=\sq 2 \de x$.  For all $w_1>0$, $\mu_V(e^{w_1|Z_1|^2})<\infty$. For $0< w_2<\ff 1 {2\de^2}$, $\mu_V(e^{w_2|Z_2|^2})<\infty$. Let $\bar Q=\left(\begin{array}{cc} -(1-\th)a  & (1+\th)a\\ (1+\th)b  &-(1-\th)b \end{array} \right)$. Then
$$-(K+\bar Q)=\left(\begin{array}{cc} 1+(1-\th)a-\ff 1 {2w_1} & -(1+\th)a\\ -(1+\th)b & 1+(1-\th)b-\ff 1 {2w_2}\end{array}\right).$$ 
By \cite[Proposition 2.4]{Shao}, $-(K+\bar Q)$ is a nonsingular M-Matrix if and only if 
$$\left\{\begin{array}{l}
1+(1-\th)a-\ff 1 {2w_1}>0\\
\left(1+(1-\th)a-\ff 1 {2w_1}\right)\left(1+(1-\th)b-\ff 1 {2w_2}\right)-(1+\th)^2ab>0.
\end{array}
\right.
$$
Consequently, if $\de$ and $\th$ satisfy 
$$\de^2<\ff {1+(1-\th)(a+b)-4\th ab} {1+(1-\th)a},$$
then there are $w_1$ and $w_2$ such that \eqref{ewZ} and \eqref{inequ_MM} hold. Moreover, if 
$$\th <\ff {1+2b+a} {a+2b+8ab},$$
then there are $\de$, $w_1$ and $w_2$ such that $\sq 2\de-1\geq 0$ and \eqref{ewZ} and \eqref{inequ_MM} hold.  That means for the state $i=2$, the diffusion of this state is not recurrent,  then it does not have an invariant probability measure. 
\end{exa}

If $P_t$ has an invariant probability measure $\mu$, then $P_t$ can be extended to be a semigroup in $L^1(\R^d\times\SS,\mu)$. Consequently, the weakly coupled elliptic system $(L,C_c^2(\R^d\times\SS))$ has one  extension in $L^1(\R^d\times\SS,\mu)$ which generates  a semigroup. For the uniqueness of the extension  in $L^1(\R^d\times\SS,\mu)$, we present the following result, where $\mu$ is introduced in Theorem \ref{thm_eu}
\beg{thm}\label{un_semigroup}
Under the assumption of Theorem \ref{thm_eu}, if in addition there exists $\ep>0$ such that
$$\max_{1\leq i\leq N}\mu_V(e^{\ep||\si_i||^2})<\infty,$$
then there is only one extension of $(L,C_c^2(\R^d\times\SS))$ that  generates a $C_0$-semigroup in $L^1(\R^d\times\SS,\mu)$.
\end{thm}

The rest  of the paper is organized  as follows. In Section 2, we shall study a elliptic system corresponding to invariant measures, called the weak coupled system of measures of \eqref{main_equ}, where some local a priori estimates of  measures will be presented. In Section 3, the entropy estimate of the density of invariant probability measures will be concern about.  In the last section, we shall prove our main results.

To make the article more concise, we shall introduce the following notations in the rest parts of the paper
\beg{enumerate}[(N1)]
\item For all $q\in [1,\infty]$, $q^*:=\1_{[1<q<\infty]}\ff {q} {q-1}+\1_{[q=1]}\infty+\1_{[q=\infty]} 0$.
\item For a Borel measurable function $f$ on $\R^d\times\SS$, despite $f(x,k)$, we also denote by $f_k(x)$ the value of $f$ at $(x,k)$. $|f|(x):=\left(\sum_{k=1}^N f_k^2(x)\right)^{1/2}$.

\item  A measure $\mu$ on $\R^d\times\SS$, we define $\mu_k(\d x)= \mu(\d x,\{k\})$.

\end{enumerate}

\section{Regularity of weak coupled system of measures}

In this section, we allow $N=\infty$. Let $\mu$ be an invariant probability measure of $P_t$ on $\sB(\R^d\times \SS)$.   Then  for  a smooth function $f$ on $\R^d\times\SS$ with compact support(there is $M\geq 0$ such that $f(x,i)=0$ if $|x|\geq M$ and $i\geq M$), It\^o's formula will yield that 
\beg{align}
P_t f(x,k)=f(x,k) +\int_0^t P_sL f(x,k)\d s.
\end{align}
Since $\mu$ is an invariant measure, 
$$\mu(P_tf)=\mu(f)+\int_0^t\mu(P_sL f)\d s=\mu(f)+\int_0^t\mu(Lf)\d s.$$
Then
$$0=\int_0^t\mu(L f)\d s,~t>0$$
and so $\mu(L f)=0$. Hence, to study the invariant measure of $P_t$, we shall start form the equation $\mu(L f)=0$. Since $\mu$ is a measure on $\R^d\times\SS$, we can view $\mu$ as a system of measures $(\mu_i)_{i\in\SS}$, where $\mu_i$ is defined as in (N3). Thus the equation $\mu(Lf)=0$ is a system of elliptic equations for measures, called weakly coupled elliptic system of measures. 

Elliptic and parabolic equations for measures have been intensively studied, and  results on equations for measures are important to the study of invariant probability measures, see  \cite{BKR,BKR2}.   Here, we shall extend some results in \cite{BKR} to the case of weakly coupled elliptic system of measures.

\beg{lem}\label{ex_den}
Assume that (H\ref{H1}) holds.  Let $p\in (d,\infty)\cap [2,\infty)$ and 
$\mu$ be a probability measure on $\R^d\times\SS$ such that
$$\mu(L f)=0,~f\in C_c(\R^d\times\SS).$$ 
Suppose in addition that for all $j\in\SS$, $Z_j\in L^p_{loc}(\R^d)\cap L_{loc}(\R^d,\d \mu_j)$, and $q_j(x)$ is locally bounded and there exists $\de(j)\in\N$ such that
\beg{align}
q_{kj}(x)=0,~|k-j|>\de(j).
\end{align} 
Then  for each $j\in\SS$, there is a continuous function $\hat\rh_j\in W_{loc}^{1,p}(\R^d)$ such that $\mu_j(\d x)=\hat\rh_j(x)\d x$.
\end{lem}
\beg{proof}
For $i\in\SS$, let $f(x,k)=f(x)\1_{[k=i]}+0\1_{[k\neq i]}$, and 
$$\nu(\d x)=-\sum_{k\neq i} q_{ki}(x)\mu_k(\d x). $$
Then
\beg{align}\label{equ_Li}
\int_{\R^d} \left(L_i^Zf(x)-q_i(x)f(x)\right)\mu_i(\d x)=\int_{\R^d} f(x)\nu(\d x).
\end{align}
According to \cite[Corollary 2.3]{BKR}, $\mu_i$ has a density $\hat\rh_i\in L^{d^*}_{loc}(\R^d)$ and $||\rh_i||_{L^{d^*}(\cO)}$ is bounded by a constant depending on $||Z_i||_{L(\cO,\d \mu_i)}$.  Then
\beg{align}\label{equ_rhi}
\int_{\R^d} \left(L^Z_if(x)-q_i(x)f(x)\right)\mu_i(\d x)=\int_{\R^d} f(x) \left(-\sum_{k\neq i}q_{ki}(x)\hat\rh_k(x)\right)\d x.
\end{align}

Let 
$$\tld Z_k^j(x)=Z_k^{0,j}(x)+\sum_{l=1}^d\si_k^{jl}(x)Z_k^l(x),~x\in\R^d,~k\in\SS,~1\leq j\leq d.$$
For all $f\in C^2_c(\R^d\times \SS)$. By $\mu(Lf)=0$ and the definition of $\hat\rh_k$, we obtain that
\beg{align}\label{equ_pde}
\sum_{k=1}^N\int_{\R^d}\left[\sum_{i,j=1}^d a_k^{ij}(x)\pp_i\pp_j f_k(x) +\sum_{j=1}^d\tld Z^j_k\pp_j f_k(x)+\sum_{j=1}^Nq_{kj}(x)f_j(x)\right]\hat\rh_k(x)&\d x\nonumber\\
&=0.
\end{align}
Fix $r\in (p^*,d^*]$. For all $q \geq  \ff {r^*p} {p-r^*}$, we have
$$\ff {q^*} p+ \ff {q^*} {r}\leq 1,\qquad q^*< r,$$
which implies that $(|\tld Z_k|+1)\hat\rh_k\in L^{q^*}_{loc}(\R^d)$  for all $k$ if $\hat\rh_k\in L^r_{loc}(\R^d)$.  Since $f\in C_c(\R^d\times \SS)$, there are $m\in\N$ such that $f_k=0$, $k>m$ and a bounded domain $\cO\subset \R^d$ such that $\bigcup_{k\leq m}\supp f_k\subset \cO$. Let 
$$\be_m(x)=\sum_{j\leq m}\max_{|k- j|\leq \de(j)}|q_{kj}|^2(x),\qquad M=\max_{j\leq m}(j+\de(j)).$$
Then \eqref{equ_pde} implies that
\beg{align}\label{grad_Sum_N}
&\left|\sum_{k=1}^N\int_{\cO}\left[\sum_{i,j=1}^d a_k^{ij}(x)\pp_i\pp_j f_k(x)\hat\rh_k(x)\right]\d x\right|\nonumber\\
&\qquad\leq \sum_{k=1}^N\int_{\cO} |\tld Z_k||\nn f_k|\hat\rh_k(x)\d x+\int_{\cO} |f|(x)\be^{\ff 1 2}_m(x)\left(\sum_{k=1}^M\hat\rh_k(x)\right)\d x\nonumber\\
&\qquad \leq \left(\sum_{k=1}^M\int_{\cO}(|\tld Z_k|\hat\rh_k)^{q^*}(x)\d x\right)^{\ff 1 {q^*}}\left(\sum_{k=1}^N\int_{\cO}|\nn f_k|^{q}(x)\d x\right)^{\ff 1 {q}}\\
&\qquad\qquad+\left[\int_{\cO} \be_m^{\ff {q^*} 2}(x)\left(\sum_{k=1}^M\hat\rh_k(x)\right)^{q^*}\d x\right]^{\ff 1 {q^*}}\left(\int_{\cO}|f|^q(x)\d x\right)^{\ff 1 {q}}\nonumber\\
&\qquad \leq C\left(\sum_{k=1}^N\int_{\cO}|\nn f_k|^{q}(x)\d x\right)^{\ff 1 {q}},\nonumber
\end{align}
where the last inequality we have used the Poincar\'e inequality(see \cite[Theorem 6.30]{AdF}),  and the constant $C$ depends on $\cO$, $\be_m$, $M$, $||\tld Z_k||_{L^p(\cO)}$, $||\hat\rh_k||_{L^{r}(\cO)}$, $\la_k$, $k=1,\cdots, M$. Fix some $k_0$. Let $f_{k}=f_{k_0}$ if $k=k_0$ and $f_{k}=0$ otherwise. Then
\beg{align}\label{grad_sum}
\left|\int_{\cO}\left[\sum_{i,j=1}^d a_{k_0}^{ij}(x)\pp_i\pp_j f_{k_0}(x)\hat\rh_{k_0}(x)\right]\d x\right|\leq C\left(\int_{\cO}|\nn f_{k_0}|^{q}(x)\d x\right)^{\ff 1 {q}}.
\end{align}

Next we can adapt a procedure introduced in \cite{BKR}. Let $x_0\in \R^d$, $R>0$, and $B_R(x_0)$ be the open ball with center  $x_0$ and radius $R$. Let $\et\in C_c^\infty(B_R(x_0))$. Then
\beg{align}\label{grad_B_R}
&\left|\int_{B_R(x_0)}\left[\sum_{i,j=1}^d a_{k_0}^{ij}(x)\pp_i\pp_j f_{k_0}(x)(\et\hat\rh_{k_0})(x)\right]\d x\right|\nonumber\\
&\qquad\leq \left|\int_{B_R(x_0)}\left[\sum_{i,j=1}^d a_{k_0}^{ij}(x)\pp_i\pp_j (\et f_{k_0})(x)\hat\rh_{k_0}(x)\right]\d x\right|\nonumber\\
&\quad\qquad+\left|\int_{B_R(x_0)}\left[\sum_{i,j=1}^d a_{k_0}^{ij}(x)\pp_i\pp_j \et(x) (f_{k_0}(x)\hat\rh_{k_0})(x)\right]\d x\right|\\
&\quad\qquad+2\left|\int_{B_R(x_0)}\left[\sum_{i,j=1}^d a_{k_0}^{ij}(x)\pp_i\et \pp_j f_{k_0}(x)\hat\rh_{k_0}(x)\right]\d x\right|\nonumber\\
&\qquad \leq C\left(\int_{B_R(x_0)}|\nn (\et f_{k_0})|^{q}(x)+|f_{k_0}|^q(x)\d x\right)^{\ff 1 {q}}\nonumber\\
&\qquad \leq C\left(\int_{B_R(x_0)}|\nn  f_{k_0}|^{q}(x)\d x\right)^{\ff 1 {q}}.\nonumber
\end{align}
the constant $C$ depends on $B_R$, $\be_m$, $M$, $\et$, $||\tld Z_k||_{L^p(B_R)}$, $||\hat\rh_k||_{L^{r}(B_R)}$, $\la_k$, $k=1,\cdots, M$. By \eqref{grad_B_R} and  \cite[Theorem 2.7]{BKR}, for $R$ small enough(independent of $x_0$), $\et\hat\rh_{k_0}\in W^{1,q^*}_{0}(B_R)$ with $1\leq q^*\leq \ff {pr} {p+r}$. The point $x_0$ is arbitrary. According to Sobolev embedding theorem, if $\ff {pr} {p+r}\geq d$, then $\hat\rh_{k_0}$ is local bounded, and if $\ff {pr} {p+r}< d$, then $\hat\rh_{k_0}\in L^{r_1}_{loc}(\R^d)$ with 
$$r_1=\ff {prd} {d-\ff {pr} {p+r}}=\ff {prd} {pd-r(p-d)}$$
which implies that
$$\ff {r_1} {r} =\ff {pd} {pd-r(p-d)}> \ff {pd} {pd-p^*(p-d)}>1.$$
Since $p>d$ yields that $p^*<d^*$,   starting from $r\in (p^*,d^*]$,  we can get  a sequence $\{r_n\}$ by this procedure such that $\hat\rh_{k_0}\in W^{1,\ff { p r_n} {p+r_n}}_{loc}(\R^d)$  until $\ff {pr_n} {p+r_n}\geq d$. Indeed, we can obtain that $\hat\rh_{k_0}\in W_{loc}^{1, l}(\R^d)$ with some $l>d$. Then  $\hat\rh_{k_0}$ is continuous and local bounded due to Sobolev embedding theorem.  From \eqref{grad_Sum_N}, \eqref{grad_sum}, $|\tld Z_k|\in L^p_{loc}(\R^d)$ and that $k_0$ is arbitrary, we obtain   for $f\in C_c(\R^d\times\SS)$
$$\left|\sum_{k=1}^N\int_{\cO}\left[\sum_{i,j=1}^d a_k^{ij}(x)\pp_i\pp_j f_k(x)\hat\rh_k(x)\right]\d x\right|\leq C\left(\sum_{k=1}^N\int_{\cO}|\nn f_k|^{p^*}(x)\d x\right)^{\ff 1 {p^*}}$$
with constant $C$ depends on $\cO$, $\be_m$, $M$, $||\tld Z_k||_{L^p(\cO)}$, $||\hat\rh_k||_{L^{r}(\cO)}$, $\la_k$, $k=1,\cdots, M$. According to \cite[Theorem 2.7]{BKR}, we also obtain that $\hat\rh_k\in W^{1,p}_{loc}(\R^d)$, and $||\hat\rh_k||_{H^{1,2}(\cO)}$ bounded by a constant depends on $p$, $d$, $\cO$, $||\tld Z_k||_{L^p(\cO)}$, $||\hat\rh_k||_{L^{d^*}(\cO)}$, $\la_k$, $a_k$, $M$, $\be_m$.

\end{proof}

Next,  we shall prove the following weak Harnack inequality, which is crucial to prove that  $$\rh_k:=\ff {e^{-V}} {\pi_k}\hat\rh_k$$
is positive. The proof mainly follows the line of \cite[Theorem 4.15]{HanL} indeed.

\beg{lem}\label{wHar}
Under the assumption of Lemma \ref{ex_den}, fixing $k\in\SS$, 
for all $0<r<R<\infty$, it  holds that
$$\mathrm{ess}\inf_{B_r} \rh_k\geq C\left(\int_{B_R}\rh_k^p\d x\right)^{\ff 1 p},~0<p<\ff d {d-2} ,$$
where $C$ depends on $\si_k$, $e^{V}$, $||q_k||_{L^{\ff p 2}(B_{R})}$, $||Z_k||_{L^{p}(B_{R})}$, $d$, $p$, $r $, $R$.
\end{lem}
\beg{proof}
By Lemma \ref{ex_den}, $\rh_k\in W_{loc}^{1,p}(\R^d)$ is  locally bounded.  \eqref{equ_rhi} and $\hat\rh_k=\pi_k\rh_k e^{V}$ imply that  
$$\int_{\R^d}(L_k^Z-q_k(x))f(x)\rh_k e^{V(x)}\d x=-\int_{\R^d}\sum_{j\neq k}q_{jk}(x)\rh_j(x) f(x) e^{V}\d x,~f\in C_c(\R^d).$$
Then by the integration by part formula, we have 
\beg{align}\label{ineq_rhk}
&\int_{\R^d}\<\si_k^*\nn f(x),\si_k^* \nn \rh_k(x)\>e^{V(x)}\d x+\int_{\R^d}q_k(x)f(x)\rh_k(x)e^{V(x)}\d x\nonumber\\
&= \int_{\R^d} \<\si_kZ_k,\nn f(x)\>\rh_k(x) e^{V}\d x+\int_{\R^d}f(x)\sum_{j\neq k}q_{jk}(x)\rh_j(x)e^{V}\d x\\
&\geq  \int_{\R^d} \<\si_kZ_k,\nn f(x)\>\rh_k(x) e^{V}\d x,~f\geq 0,~f \in C_c(\R^d).\nonumber
\end{align}
The remainder of the proof will start from this inequality.

Let $\de>0,~\be>0$ and $f=\ps^2(\rh_k+\de)^{-(\be+1)}$ with $\ps\in C_c(\R^d)$.  Then
\beg{align*}
&2\int_{\R^d}\ff {\<\si_k^*\nn\ps,\si^*_k\nn \rh_k\>\ps} {(\rh_k+\de)^{\be+1}} e^{V(x)}\d x-(\be+1)\int_{\R^d}\ff {|\si_k^*\nn \rh_k|^2\ps^2} {(\rh_k+\de)^{\be+2}} e^{V(x)}\d x\\
&\geq -\int_{\R^d}\ff {q_k(x) \ps^2(x) \rh_k} {(\rh_k+\de)^{\be+1}} e^{V(x)}\d x+2\int_{\R^d}\ff {\<\si_kZ_k,\nn\ps\>\ps\rh_k} {(\rh_k+\de)^{\be+1}} e^{V(x)}\d x\\
&\qquad-(\be+1)\int_{\R^d}\ff {\<\si_kZ_k,\nn\rh_k\>\ps^2\rh_k} {(\rh_k+\de)^{\be+1}}e^{V(x)}\d x.
\end{align*}    
Thus
\beg{align*}
&\ff {4(\be+1)} {\be^2}\int_{\R^d}|\si_k^*\nn(\rh_k+\de)^{-\ff {\be} 2}|^2 \ps^2 e^{V(x)}\d x\\
&\leq \ff {4} {\be} \int_{\R^d}\<\si_k^*\nn\ps,\si_k^*\nn(\rh_k+\de)^{-\ff {\be} 2}\> \ps(\rh_k+\de)^{-\ff {\be} 2} e^{V(x)}\d x+\int_{\R^d}\ff {q_k(x) \ps^2} {(\rh_k+\de)^{\be}} e^{V(x)}\d x\\
&\quad +2\int_{\R^d}\ff {|Z_k||\si_k^*\nn\ps||\ps|} {(\rh_k+\de)^{\be}} e^{V(x)}\d x\\
&\quad+\ff {2(\be+1)} {\be}\int_{\R^d} \ff {|Z_k|\ps^2|\si_k^*\nn(\rh_k+\de)^{-\ff {\be} 2}|} {(\rh_k+\de)^{\ff {\be} 2 } } e^{V(x)}\d x.
\end{align*}
Let $w=(\rh_k+\de)^{-\ff {\be} 2}$. Then  from the  inequality above, we  can obtain that 
 \beg{align*}
&\int_{\R^d}|\si_k^*\nn w|^2\ps^2 e^{V(x)}\d x\\
&\leq \ff {\be} {4(\be+1)}\int_{\R^d} |\si^*_k\nn \ps||\si_k^*\nn w|\ps w e^{V(x)}\d x+\ff {\be^2} {4(\be+1)}\int_{\R^d} q_k(x)\ps^2 w^2 e^{V}\d x\\
&+\ff {\be^2} {2(\be+1)}\int_{\R^d}|Z_k||\si_k^*\nn\ps|\ps w^2 e^{V(x)}\d x+\ff {\be} {2} \int_{\R^d} | Z_k| |\si_k^* \nn w| \ps^2 we^{V(x)}\d x.
 \end{align*}
Thus
\beg{align*}
\int_{\R^d} |\si_k^*\nn w|^2\ps^2 e^{V}\d x&\leq (1+\be)\int_{\R^d} |\si_k^*\nn\ps|^2 w^2 e^{V}\d x+2\be\int_{\R^d}q_k(x)\ps^2 w^2 e^{V(x)}\d x\\
&\quad +\be\int_{\R^d} |Z_k|^2\ps^2 w^2 e^{V}\d x.
\end{align*}
Suppose $\supp \ps\subset \cO$. Letting $\mathrm{osc}_{\cO} e^{V}= \ff {\sup_{\cO} e^V} {\inf_{\cO} e^V}$ and $\La_\cO=\sup_{\cO}||\si^*_k||^2$, we obtain that
\beg{align}\label{ineq_wps}
&\int_{\R^d} |\nn (w\ps)|^2 \d x\nonumber\\
&\leq C(\la_k,\mathrm{osc}_{\cO} e^{V},\La_{\cO})(1+\be)\int_{\R^d}\left( |\nn\ps|^2+q_k(x)\ps^2+|Z_k|^2\ps^2\right) w^2\d x.
\end{align} 
By H\"older inequality,  
\beg{align*}
\int_{\R^d} (q_k(x)+|Z_k|^2)\ps^2 w^2 \d x\leq \left(\int_{\cO} (q_k(x)+|Z_k|^2)^{\ff p 2} \d x\right)^{\ff 2 p}\left(\int_{\R^d} (\ps w)^{\ff {2p} {p-2}}\d x\right)^{\ff {p-2} {p}}.
\end{align*}
Since $p>d$, H\"older inequality and Sobolev inequality imply that  for all $\ep>0$
\beg{align*}
||\ps w||_{L^{\ff {p} {p-2}}(\cO)}&\leq \ep||\ps w||_{L^{\ff d {d-1}}(\cO)}+C(d, p)\ep^{-\ff d {p-d}}||\ps w||_{L^2(\cO)}\\
& \leq \ep C(d, \cO)||\nn (\ps w)||_{L^2(\cO)}+C(d,p)\ep^{-\ff d {p-d}} ||\ps w||_{L^2(\cO)}.
\end{align*}
Substituting this into \eqref{ineq_wps} and choosing some small $\ep$,  we obtain that 
\beg{align*}
&\int_{\R^d} |\nn (w\ps)|^2 \d x\leq C(1+\be)^{\al}\int_{\R^d}\left( |\nn\ps|^2+\ps^2\right) w^2\d x
\end{align*}
with some $\al>0$ and $C$ depending on $\la_k,\mathrm{osc}_{\cO} e^{V},\La_{\cO}, ||q_k||_{L^{\ff p 2}(\cO)}, ||Z_k||_{L^{p}(\cO)},p,d$.
By the Sobolev embedding theorem, 
\beg{align*}
\left(\int_{\R^d}(\ps w)^{2\ga}\d x\right)^{\ff 1 {\ga}} \leq C(1+\be)^\al \int_{\R^d} \left(|\nn\ps|^2 + \ps^2\right) w^2\d x
\end{align*}
where $\ga=\ff {d} {d-2}$ for $d\geq 3$ and $\ga>2$ for $d= 2,~1$. Choosing $\ps$ such that  $\1_{B_{r_1}}\leq \ps \leq \1_{ B_{r_2}}$ and $|\nn\ps|\leq \ff 2 {r_2-r_1}$, where $B_{r_1}$ and $B_{r_2}$ are some balls with radius $r_1$ and $r_2$, which are subsets of $\cO$. We arrive at 
\beg{align}
\left(\int_{B_{r_1}}w^{2\ga } \d x\right)^{\ff 1 {\ga}}\leq \ff {C(1+\be)^\al} {(r_2-r_1)^2} \int_{B_{r_2}} w^2\d x.
\end{align}
Thus
\beg{align}
\left(\int_{B_{r_1}} (\rh_k+\de)^{-\be\ga}\d x\right)^{\ff {1 } {\be\ga}}\leq \left(\ff {C(1+\be)^\al} {(r_2-r_1)^2}\right)^{\ff 1 {\be}} \left(\int_{B_{r_2}} (\rh_k+\de)^{-\be}\d x\right)^{\ff 1 {\be}}.
\end{align}
Hence, we can get the inverse H\"older inequality for $(\rh_k+\de)^{-1}$ by iteration(see \cite{HanL} for instance)
\beg{align}\label{inver_H1}
||(\rh_k+\de)^{-1}||_{L^{p_1}(B_{r_1})}\leq C||(\rh_k+\de)^{-1}||_{L^{p_2}(B_{r_2})}
\end{align}
with any $p_1>p_2>0$ and $B_r\subset B_R$, where the constant $C$ depends on $\la_k$, $\mathrm{osc}_{B_{r_2}} e^{V}$, $\La_{B_{r_2}}$, $||q_k||_{L^{\ff p 2}(B_{r_2})}, ||Z_k||_{L^{p}(B_{r_2})},d,p,r_1,r_2$. Moreover, letting $p_1\ra \infty$, we obtain 
\beg{align}\label{inf}
\mathrm{ess}\inf_{B_{r_1}}\rh_k +\de\geq C ||(\rh_k+\de)^{-1}||^{-1}_{L^{p_2}(B_{r_2})}.
\end{align}

Next, we shall prove that for all $0<q<\ff d {d-2}$ and $R>0$, there is $C>0$ such that
\beg{align}\label{q-q}
\int_{B_R}(\rh_k+\de)^{-q}\d x\int_{B_R}(\rh_k+\de)^{q}\d x\leq C.
\end{align}
Let $f=(\rh_k+\de)^{-1}\ps^2$ in \eqref{ineq_rhk} and $w=\log(rh_k+\de)$.  Then
\beg{align*}
\int_{\R^d}|\si_k^*\nn w|^2\ps^2 e^{V}\d x&\leq  \int_{\R^d} |\si_k^*\nn\ps||\si_k^*\nn\ps|\ps e^{V}\d x+\int_{\R^d}|Z_k||\nn w|\ps^2 e^{V}\d x\\
&\quad +\int_{\R^d} |Z_k||\nn\ps|\ps e^V\d x+\int_{\R^d}q_k(x)\ps^2\d x.
\end{align*}
Hence
\beg{align*}
\int_{\R^d}|\si_k^*\nn w|^2\ps^2 e^{V}\d x\leq 4 \int_{\R^d} |\si_k^*\nn\ps|^2 e^{V}\d x+4\int_{\R^d}(|Z_k|^2+q_k(x))\ps^2 e^{V}\d x.
\end{align*} 
Choosing $\ps$ such that $\1_{B_r}\leq \ps\leq  \1_{B_{2r}}\leq \1_{\cO}$ and $|\nn\ps|\leq \ff 2 r$, we obtain that there is a constant $C>0$ depending on $\la_k$, $\La_{\cO}$, $\cO$  such that
$$\int_{B_r}|\nn w|^2\d x\leq C\left( r^{d-2}+\int_{B_{2r}}\left(|Z_k|^2+q_k(x)\right)\d x\right).$$
Since 
\beg{align*}
\int_{B_{2r}}\left(|Z_k|^2+q_k(x)\right)\d x&\leq \left|\left||Z_k|^2+q_k(x)\right|\right|_{L^{\ff d 2}(B_{2r})}\left(\int_{B_{2r}}\d x\right)^{\ff {d-2} d}\\
&= 2^{d-2}r^{d-2}\left|\left||Z_k|^2+q_k(x)\right|\right|_{L^{\ff d 2}(B_{2r})},
\end{align*}
We have that 
$$\int_{B_r}|\nn w|^2\d x\leq Cr^{d-2}$$
with $C$ depends on $\la_k$, $\La_{\cO}$, $\cO$, $||q_k||_{L^{\ff d 2}(B_{2r})}$, $||Z_k||_{L^{d}(B_{2r})}$, $d$. As in \cite[Theorem 4.14]{HanL}, the Poincar\'e inequality on $B_{2r}$ and John-Nirenberg Lemma imply that there exist $p_0>0$ and $C>0$ depending on $d$, $r$ such that
$$\int_{B_r}e^{p_0|w-\bar w|}\d x\leq C$$
where $\bar w=\ff 1 {|B_r|}\int_{B_r}w\d x$. Hence
\beg{align}\label{p0-p0}
\int_{B_r} (\rh_k+\de)^{-p_0}\d x \int_{B_r} (\rh_k+\de)^{p_0}\d x&=\int_{B_r} e^{-p_0w}\d x \int_{B_r} e^{p_0w}\d x\nonumber\\
&\leq \left(\int_{B_r} e^{p_0|w-\bar w|}\d x\right)^2\\
&\leq C^2.\nonumber
\end{align}
To prove \eqref{q-q}, let $f=\ps^2(\rh_k+\de)^{-\be}$, $\be\in(0,1)$. Then, as above, we can prove the following inverse H\"older inequality
$$||\rh_k+\de||_{L^{p_1}(B_{r_1})}\leq C||\rh_k+\de||_{L^{p_2}(B_{r_2})},~0<r_1<r_2,~0<p_2<p_1<\ff d {d-2}$$
where the constant $C$ depends on $\la_k$, $\mathrm{osc}_{B_{r_2}} e^{V}$, $\La_{B_{r_2}}$, $||q_k||_{L^{\ff p 2}(B_{r_2})}$, $||Z_k||_{L^{p}(B_{r_2})}$, $d$, $p$, $r_1$, $r_2$.  Combing this with \eqref{p0-p0} and \eqref{inver_H1}, \eqref{q-q} is proved.

Lastly, \eqref{inf} and \eqref{q-q} yield that
\beg{align*}
\mathrm{ess}\inf_{B_{r_1}}\rh_k +\de &\geq C ||(\rh_k+\de)^{-1}||^{-1}_{L^{p_2}(B_{r_2})}\\
&\geq  C\left(||(\rh_k+\de)^{-1}||_{L^{q}(B_{r_2})}||(\rh_k+\de)||_{L^{q}(B_{r_2})}\right)^{-1} ||(\rh_k+\de)||_{L^{q}(B_{r_2})}\\
&\geq C||(\rh_k+\de)||_{L^{q}(B_{r_2})}.
\end{align*}
Letting $\de\ra 0^+$,  we obtain the weak Harnack inequality.

\end{proof}
 
\beg{rem}\label{rem-add}
If $Q$ is fully coupled and $\mu$ is a probability measure, then it follows from the weakly Harnack inequality that $\rh_k>0$ for all $k\in\SS$. Indeed, the weakly Harnack inequality implies that if $\mu_k$ is a positive measure, then $\rh_k$ is positive on $\R^d$. Since $\mu$ is a probability measure, there exist $k\in\SS$ such that $\rh_k\pi_ke^{V(x)}\d x$ is a positive measure. Let $\SS_1$ be the set of all $\mu_j$ such that $\mu_j$ is zero, and $\SS_2=\SS-\SS_1$. If $\SS_1$ is non-empty, then it follows from \eqref{equ_rhi} that 
$$q_{l_1l_2}=0, \mbox{a.e.},\qquad l_1\in\SS_1,l_2\in\SS_2.$$
It  contradicts that $Q$ is fully coupled on $\R^d$.
\end{rem}

\section{Entropy estimates of $\rh$}
In this section, we shall study the the entropy estimate of the density $\rh$. Firstly, we shall prove a Girsanov's theorem, which will be used to prove the existence of the invariant measure.
Let $(X_t,\Lambda_t)$ be the solution of the following equation
\beg{align}
\d X_t&=b^{(1)}_{\La_t}(X_t)\d t+\sq 2g_{\La_t}(X_t)\d W_t,~X_0=\xi \label{BasicEqu1}\\
\d \La_t&=\int_0^\infty h(X_t,\La_{t-},z)N(\d z,\d t),~\La_0=\la.\label{BasicEqu2}
\end{align}
Let 
\beg{align}
\tld W_t&=W_t-\ff 1 {\sq 2}\int_0^t b^{(2)}_{\La_s}(X_s)\d s\\
R_t&=\exp\Big\{\ff 1 {\sq 2}\int_0^t\Big\<b^{(2)}_{\La_s}(X_s),\d W_s\Big\>-\ff 1 2\int_0^t|b^{(2)}_{\La_s}(X_s)|^2\d s\Big\}.
\end{align}
Let $\{\mathbf{p}_t\}_{t\geq 0}$ be the Poisson point process corresponding to the random measure $N(\d z,\d t)$, i.e. $\mathbf{p}$ is a Poisson point process which takes  value in $(0,\infty)$ and is independent of the Brownian motion $\{W_t\}_{t\geq 0}$ such that 
$$N(U,(0,t])(\om)=\#\{s\in D_{\mathbf{p}(\om)}~|~s\leq t,~\mathbf{p}_s(\om)\in U\},$$
where $D_{\mathbf{p}(\om)}$ is the domain of the point function $\mathbf{p}(\om)$.

\beg{lem}\label{GirTh}
Suppose that \eqref{BasicEqu1}-\eqref{BasicEqu2} has a pathwise unique non-explosive solution, and $\{R_t\}_{t\in[0,T]}$ is a martingale and 
$$\P\left(\int_0^T| g_{\La_s}(X_s) b^{(2)}_{\La_s}(X_s)|^2\d s<\infty\right)=1.$$ 
Then, under $R_T\P$,  $\{\tld W_t\}_{t\in[0,T]}$  is a Brownian process and $\{\mathbf{p}_t\}_{t\geq 0}$ is a Poisson point process such that  they are mutually independent. Consequently,  the following equation has a weak solution 
\beg{align*}
\d X_t&=b^{(1)}_{\La_t}(X_t)\d t+\sq 2g_{\La_t}(X_t)\d W_t+\sq 2 g_{\La_t}(X_t) b^{(2)}_{\La_t}(X_t)\d t,~X_0=\xi  \\
\d \La_t&=\int_0^\infty h(X_t,\La_{t-},z)N(\d z,\d t),~\La_0=\la.
\end{align*}
\end{lem}
\beg{proof} 
$\{R_t\}_{t\in[0,T]}$ is a martingale, $\tld W_t$ is a Brownian motion due to Girsanov's theorem. According to \cite[Theorem I.6.3]{IkW}, we only need to prove  that the compensator of $N(\d z,\d t)$ under $R_T\P$ is $\d z\d t$. 

Let  $K$ be a measurable subset of $\R^+$ with its Lebesgue measure $|K|<\infty$, $N(K,t):=N(K,(0,t]),~t\in[0,T]$. Then
\beg{align*}
&\E\Big[R_T\left(N(K,t)-|K|t\right)\Big|\sF_s\Big]\\
&=-\E\left[|K|tR_T\Big|\sF_s\right]+\E\left[N(K,s)R_T\Big|\sF_s\right]+\E\left[\left(N(K,t)-N(K,s)\right)R_T\Big|\sF_s\right]\\
&=\left(N(K,s)-|K|t\right)R_s+\E\left[\left(N(K,t)-N(K,s)\right)R_t\Big|\sF_s\right]\\
&=\left(N(K,s)-|K|t\right)R_s+R_s\E\left[\left(N(K,t)-N(K,s)\right)R_{s,t}\Big|\sF_s\right],
\end{align*}
where
\beg{align*}
R_{s,t}&=\exp\Big\{\int_s^t\Big\< b^{(2)}_{\La_r}(X_r),\d W_r\Big\>-\ff 1 {\sq 2} \int_s^t| b^{(2)}_{\La_r}(X_r)|^2\d r\Big\}.
\end{align*}
Let $\sG=\sF_s\vee \si\Big(N(U,r)~\Big|~r\leq t,~U\in\sB(\R^+)\Big)$. Then 
\beg{align*}
\E\Big[\left(N(K,t)-N(K,s)\right)R_{s,t}\Big|\sF_s\Big]=\E\Big[\left(N(K,t)-N(K,s)\right)\E\Big[R_{s,t}|\sG\Big]\sF_s\Big].
\end{align*}
Let  $\Pi_{[s,t]}$ be the totality of point functions on $[s,t[$ for $s\leq t$, and let $\sD_{s,t}(\SS)$ be   all the cadlag function on $[s,t]$ with value in $\SS$. Since the equation has a unique non-explosive strong solution, there exists $F:\R^d\times \SS\times C([s,t],\R^d)\times \Pi_{[s,t]}\ra C([s,t],\R^d)\times \sD_{s,t}(\SS)$ such that 
$$(X_r,\La_r)=F(X_s,\La_s,W_{[s,t]},\mathbf {p}_{[s,t]})_r,~r\in [s,t].$$
Thus
\beg{align*}
&\E\Big[R_{s,t}\Big|\sG\Big]\\
&=\E\Big[\exp\Big\{\int_s^t\Big\<b^{(2)}(F(X_s,\La_s,W_{[s,t]},\mathbf {p}_{[s,t]})_r),\d W_r\Big\>\\
&\quad -\ff 1 {\sq 2}\int_s^t| b^{(2)}(F(X_s,\La_s,W_{[s,t]},\mathbf {p}_{[s,t]})_r)|^2\d r\Big\}\Big|\sG\Big]\\
&=\E\exp\Big\{\int_0^{t-s}\Big\<( b^{(2)}(F(x,k,\th_s( W_{[0,t-s]})+x,p_{[s,t]})_r),\d W_u\Big\>\\
&\quad -\ff 1 {\sq 2}\int_0^{t-s}| b^{(2)}(F(x,k,\th_s( W_{[0,t-s]})+x,p_{[s,t]})_r)|^2\d u\Big\}\Big|_{(x,k,p)=(X_s, \La_s, \mathbf{p})}\\
&\leq 1,\quad \P(\cdot|\sG)\mbox{-a.s.}
\end{align*}
Combining with $\E R_{s,t}=\E\left(\E\Big[R_{s,t}\Big|\sG\Big]\right)=1$, we have that 
$$\E\Big[R_{s,t}\Big|\sG\Big]=1,\quad \P(\cdot|\sG)\mbox{-a.s.}$$
Since $\sF_s\subset \sG$, we obtain that
\beg{align*}
\E\Big[\left(N(K,t)-N(K,s)\right)\E\Big[R_{s,t}|\sG\Big]\Big|\sF_s\Big]=\E\Big[N(K,t)-N(K,s)\Big|\sF_s\Big]=|K|(t-s).
\end{align*}
Then
$$\E\Big[R_T\Big(N(K,t)-|K|t\Big)\Big|\sF_s\Big]=R_s\Big[N(K,s)-|K|s\Big].$$
Thus for $A\in \sF_s$
\beg{align*}
&\E_{R_T}\Big\{\left(N(K,t)-|K|t\right)\1_{A}\Big\}=\E\Big\{ R_T\left(N(K,t)-|K|t\right)\1_{A}\Big\}\\
&=\E \Big\{R_s\Big[N(K,s)-|K|s\Big]\1_{A}\Big\}=\E \Big\{R_T\Big[N(K,s)-|K|s\Big]\1_{A}\Big\}\\
&=\E_{R_T} \Big\{\Big[N(K,s)-|K|s\Big]\1_{A}\Big\}.
\end{align*}
Therefore
$$\E_{R_T}\Big[N(K,t)-|K|t~\Big|\sF_s\Big]=N(K,s)-|K|s,$$
which implies that the compensator of $N(\d z,\d t)$ under $R_T\P$ is $\d z\d t$.

\end{proof}

Let $P_t^{0}$ be the Markov semigroup generated by the Dirichlet form $\sE^0$.
Then $L^0$ is  the generator of $P_t^0$, and $\mu_{\pi,V}$ is the invariant probability measure of $P_t^0$. Moreover, for $f\in L^2(\R^d\times\SS,\mu_{\pi,V})$, 
$$(P_t^0)f(x,i)=P_t^if_i(x), x\in\R^d,~i\in\SS,$$
where $P_t^{i}$ is the semigroup generated by $L^0_i$, and $P_t^0$ can be extended to a contrastive $C_0$-semigroup in $L^{p}(\R^d\times\SS,\mu_{\pi,V})$, $p\in[1,\infty]$. By (H\ref{H2}),
\beg{align*}
&\sum_{i=1}^N\pi_i\int_{\R^d} \left|\sum_{j=1}^N q_{ij}(x)f_j(x)\right|^p\mu_V(\d x)\\
&\leq \sup_{x\in\R^d}\left(\sum_{i=1}^N\pi_i\left(\sum_{j=1}^N\ff {|q_{ij}|^{\ff p {p-1}}(x)} {\pi_j^\ff {1} {p}}\right)^{p-1} \right)\int_{\R^d} \sum_{j=1}^N\pi_j f_j^p(x)\mu_V(\d x).
\end{align*}
Thus $Q$ is a bounded operator in $L^p(\R^d\times\SS,\mu_{\pi,V})$ for all $p\in[1,\infty]$.  Let  $L^Q$ be  defined as follows
$$(L^Qf)(x,i)=(L_i^0f_i)(x)+\sum_{j=1}^Nq_{ij}(x)f_j(x),~f\in C_c^2(\R^d\times\SS).$$
Then $L^Q$ generates a unique contractive $C_0$-semigroup in $L^p(\R^d\times\SS,\mu_{\pi,V})$, say $P_t^Q$.  By (H\ref{H1}), (H\ref{H2}), and  for all $f\in C^1_c(\R^d\times\SS)$ we have
\beg{align*}
&\sum_{i}\pi_i\sum_{j}q_{ij}(f_j^+\we 1)(f_i-f_i^+\we 1)\\
&=\sum_{i}\pi_i\sum_{j\neq i}q_{ij}\left[(f_j^+\we 1)f_i^-\right]-\sum_{i}\pi_i\sum_{j\neq i}q_{ij}(f_j^+\we 1)(f_i-1)^+\\
&\qquad +\sum_{i}\pi_iq_i(f_i^+\we 1)(f_i-1)^+\\
&=\sum_{i}\pi_i\sum_{j\neq i}q_{ij}\left[(f_j^+\we 1)f_i^-+(f_i^+\we 1)(f_i-1)^+\right]\\
&\qquad-\sum_{i}\pi_i\sum_{j\neq i}q_{ij}(f_j^+\we 1)(f_i-1)^+\\
&\geq \sum_{i}\pi_i\sum_{j\neq i}q_{ij}(1-f_j^+\we 1)(f_i-1)^+\geq 0.
\end{align*}
According to the proof of \cite[Theorem 3.1]{ChenZ},  $P_t^Q$ is also a  Markov semigroup. Moreover, we have 
\beg{lem}\label{hyper}
For $t>0$, $P_t^Q$ is hyperbounded, and 
$$||P_t^Q||_{L^q(\mu_{\pi,V})\ra L^{q(t)}(\mu_{\pi,V})}<\infty,~q>1,~q(t)=1+(q-1)e^{\ff {4t} {\ga}},$$ 
recalling that $\ga=\displaystyle\max_{1\leq i\leq N}\ga_i$.
\end{lem}
\beg{proof}
By \cite[Theorem 3.1.1]{Pa}, we have
\bequ\label{Pt_Pt0}
P_t^{Q}f=P_t^{0}f+\int_0^tP_{t-s}^{0}\left(QP_s^{Q}f\right)\d s,~f\in L^\infty(\R^d\times\SS,\mu_{\pi,V}).
\enqu
Then  for all $p\geq 1$
\beg{align*}
\mu_{\pi,V}(|P_t^Qf|^{p})&\leq 2^{p}\mu_{\pi,V}(|P_t^0f|^p)+t^{p-1}\int_0^t\mu_{\pi,V}(|QP_s^Qf|^{p})\d s\\
&\leq 2^{p-1}\mu_{\pi,V}(|P_t^0f|^{p})+2^{p-1} t^{p-1}\int_0^t\mu_{\pi,V}(|QP_s^Qf|^{p})\d s\\
&\leq 2^{p-1}\mu_{\pi,V}(|P_t^0f|^{p})+2^{p-1} t^{p-1}b_p\int_0^t\mu_{\pi,V}(|P_s^Qf|^{p})\d s,
\end{align*}
where $b_p=||Q||_{L^{p}}^{p}$. Gronwall's inequality yields 
$$\mu_{\pi,V}(|P_t^Qf|^{p})\leq 2^{p-1}e^{2^{p-1}t^{p}b_p}\mu_{\pi,V}(|P_t^0f|^p).$$
According to \cite[Theorem 5.1.4]{Wang04B}, (H\ref{H3}) and Remark \ref{rem_logSob} imply that $P_t^0$ is hyperbounded:
$$||P_t^{0}||_{L^{q}(\mu_{\pi,V})\ra L^{q(t)}(\mu_{\pi,V})}\leq \exp\left[\be\left(\ff 1 q-\ff 1 {q(t)}\right)\right], t>0,~q>1,$$
where $q(t):=1+(q-1)e^{\ff {4t} {\ga}}$ and $\be=\max_k(\be_k-\log {\pi_k})$. Thus 
$$\left(\mu_{\pi,V}(|P_t^Qf|^{q(t)})\right)^{\ff 1 {q(t)}}\leq 2^{1-\ff 1 {q(t)}}\exp\left(\ff {\be(q(t)-q)} {qq(t)}+\ff {2^{q(t)-1}t^{q(t)}b_{q(t)}} {q(t)}\right)\left(\mu_{\pi,V}(f^q)\right)^{\ff 1 q},$$
and the  proof of the lemma is completed.

\end{proof}

Let $(X_t,\La_t)$ be the process generated by $L^Q$, which satisfies the stochastic differential \eqref{main_equ} with $b_k(x)=Z_k^0(x)$, $k\in\SS$, $x\in\R^d$. The following lemma is a similar result of \cite[Theorem 4.1]{Wang18} for \eqref{main_equ}.

\beg{lem}\label{lem_mu}
Assume (H\ref{H1})-(H\ref{H3}). Let $p_1=\ff 1 2 (1+e^{\ff 4 {\ga}})$ and
$$c_1=\left(\ff {e^{b_1}-1} {b_1}\right)^{\ff 1 2} ||P_1^Q||_{L^2(\mu_{\pi,V})\ra L^{2p_1}(\mu_{\pi,V})},$$
where $b_1=||Q||_{L^1(\mu_{\pi,V})}$. 
There is $w>\ff {3p_1-1} {2(p_1-1)}$ such that 
\bequ\label{inequ_int}
\mu_{\pi,V}\left(e^{w|Z|^2}\right):=\sum_{i=1}^N\pi_i\mu_{V}\left(e^{w|Z_i|^2}\right)<\infty.
\enqu
Then $P_t$ has an  invariant measure $\mu=\rh\mu_{\pi,V}$ such that
$$\mu_{\pi,V}(\rh\log\rh)\leq \ff {(3p_1-1)\log\mu_{\pi,V}\left(e^{w|\si^{-1}Z|^2}\right)+4wp_1\log c_1} {2\la(p_1-1)-(3p_1-1)}.$$
\end{lem}
\beg{proof}
We first consider the Feynman-Kac semigroup for any $F\in \sB_b(\R^d\times\SS)$:
$$(P_t^Ff)(x,i)=\E\left[f(X_t^{(x,i)},\La_t^{(x,i)})\exp\left(\int_0^t F(X_s^{(x,i)},\La_s^{(x,i)})\d s\right)\right],~f\in \sB_b(\R^d\times\SS).$$
Let $p=1+\ff {p_1-1} {2p_1}$. By Lemma \ref{hyper}, $||P_1^Q||_{L^2(\mu_{\pi,V})\ra L^{2p_1}(\mu_{\pi,V})}<\infty$.
Then
\beg{align*}
&\mu_{\pi,V}(|P_1^Ff|^{2p})\leq \mu_{\pi,V}\left(\left(\E\left[\left|f(X_1,\La_1)\right|\exp\left(\int_0^1F(X_s,\La_s\d s\right)\right]\right)^{2p}\right)\\
&\leq \left(\mu_{\pi,V}\left((P_1^Qf^{p})^{2p_1}\right)\right)^{\ff 1 {p_1}}\left\{\mu_{\pi,V}\left(\left(\E e^{\ff p {p-1}\int_0^1 F(X_s,\La_s)\d s}\right)^{\ff {2(p-1)p_1} {p_1-1}}\right)\right\}^{\ff {p_1-1} {p_1}}\\
&\leq ||P_1^Q||_{L^2(\mu_{\pi,V})\ra L^{2p_1}(\mu_{\pi,V})}^2\mu_{\pi,V}(f^{2p})\left\{\int_0^1\mu_{\pi,V}\left(\E e^{\ff p {p-1}F(X_s,\La_s)}\right)\d s\right\}^{\ff {p_1-1} {p_1}}\\
&= ||P_1^Q||_{L^2(\mu_{\pi,V})\ra L^{2p_1}(\mu_{\pi,V})}^2\mu_{\pi,V}(f^{2p})\left\{\int_0^1\mu_{\pi,V}\left(P_s^Qe^{\ff p {p-1}F}\right)\d s\right\}^{\ff {p_1-1} {p_1}}\\
&\leq ||P_1^Q||_{L^2(\mu_{\pi,V})\ra L^{2p_1}(\mu_{\pi,V})}^2\mu_{\pi,V}(f^{2p})\left\{\int_0^1e^{b_1s}\mu_{\pi,V}\left(e^{\ff p {p-1}F}\right)\d s\right\}^{\ff {p_1-1} {p_1}}.
\end{align*}
Thus
$$||P_1^F||_{L^{2p}(\mu_{\pi,V})}\leq c_1^{\ff 1 p}\left(\mu_{\pi,V}\left(e^{\ff p {p-1}F}\right)\right)^{\ff {p_1-1} {2p_1p}}=c_1^{\ff 1 p}\left(\mu_{\pi,V}\left(e^{\ff {3p_1-1} {p_1-1}F}\right)\right)^{\ff {p_1-1} {3p_1-1}},$$
where $c_1=\left(\ff {e^{b_1}-1} {b_1}\right)^{\ff 1 2} ||P_1^Q||_{L^2(\mu_{\pi,V})\ra L^{2p_1}(\mu_{\pi,V})}$. By the Markov property, we have that
$$\E e^{\int_0^n F(X_s^{\mu_{\pi,V}},\La_s^{\mu_{\pi,V}})\d s}=\mu_{\pi,V}(P_n^F\1)\leq ||P_1^F||^n_{L^{2p}(\mu_{\pi,V})}\leq c_1^{\ff n p}\left(\mu_{\pi,V}\left(e^{\ff {3p_1-1} {p_1-1}F}\right)\right)^{\ff {n(p_1-1)} {2p_1p}}.$$

Next, let 
\beg{align*}
R_t^{\mu_{\pi,V}}=\exp\left\{\int_0^t\<Z(X^{\mu_{\pi,V}}_s,\La^{\mu_{\pi,V}}_s),\d W_s\>-\ff 1 2\int_0^t|Z|^2(X^{\mu_{\pi,V}}_s,\La^{\mu_{\pi,V}}_s)\d s\right\}.
\end{align*}
Then, similar to \cite[Theorem 4.1]{Wang18}, 
\beg{align*}
&\E R_n^{\mu_{\pi,V}}\log R_n^{\mu_{\pi,V}}=\ff 1 2\E \left(R_n^{\mu_{\pi,V}}\int_0^n|Z(X_s,\La_s)|^2\d s\right)\\
&\leq \ep\E R_n^{\mu_{\pi,V}}\log R_n^{\mu_{\pi,V}}+\ep \log\E\left[e^{\ff 1 {2\ep}\int_0^n|Z|^2(X_s,\La_s)}\d s\right],~\ep\in(0,1).
\end{align*}
Thus
$$\E R_n^{\mu_{\pi,V}}\log R_n^{\mu_{\pi,V}}\leq \ff {n\ep} {p(1-\ep)}\log\left[c_1\left(\mu_{\pi,V}\left(e^{\ff {3p_1-1} {2\ep(p_1-1)}|Z|^2}\right)\right)^{\ff {p_1-1} {2p_1}}\right],~\ep\in(0,1).$$
By Lemma \ref{GirTh}, 
\beg{align*}
\nu_n(f)&:=\ff 1 n\int_0^n\mu_{\pi,V}\left(P_sf\right)\d s=\ff 1 n\int_0^n \E R_n^{\mu_{\pi,V}} f(X_s^{\mu_{\pi,V}},\La_s^{\mu_{\pi,V}})\d s\\
&\leq \ff 1 n \E R_n^{\mu_{\pi,V}}\log R_n^{\mu_{\pi,V}}+\ff 1 n \log\E e^{\int_0^n f(X_s^{\mu_{\pi,V}},\La_s^{\mu_{\pi,V}})\d s}\\
&\leq \ff {\ep} {p(1-\ep)}\log\left[c_1\left(\mu_{\pi,V}\left(e^{\ff {3p_1-1} {2\ep(p_1-1)}|Z|^2}\right)\right)^{\ff {p_1-1} {2p_1}}\right]\\
&\quad+\ff 1 p\log\left[c_1\left(\mu_{\pi,V}\left(e^{\ff {3p_1-1} {p_1-1}f}\right)\right)^{\ff {p_1-1} {2p_1}} \right],~f\geq 0,~f\in\sB_b(\R^d\times\SS).
\end{align*}
By these estimates above, following the proof of \cite[Theorem 4.1]{Wang18}, we can prove the theorem. 

\end{proof}

\beg{rem}
The inequality \eqref{inequ_int} implies that for each $i\in\SS$, $Z_i$ has nice integrability w.r.t $\mu_V$ indeed.  It does not work well for \eqref{main_equ}. In the rest part of this section, we shall give a weaker condition to get a priori estimate to the  entropy of the density of invariant probability measures.
\end{rem}


We define $x\log x=0$ if $x=0$. Then we have 

\beg{lem}\label{ineq_nn_en}
Assume (H\ref{H1})-(H\ref{H3}) and that for all $i\in\SS$, $Z_i$ is bounded with compact support. 
Then there is positive $\rh\in L^1(\R^d\times\SS,\,u_{\pi,V})$ such that $\sq{\rh_i}\in H^{1,2}_{\si_i}(\R^d)$ with
\beg{align*}
&2\pi_i\mu_V(|\si^*_i\nn\sq{\rh_i}|^2)+ \pi_i\mu_V(q_i\rh_i\log\rh_i)\\
&\leq \ff {\pi_i} 2\mu_V(|Z_i|^2\rh_i)+\sum_{k\neq i}\pi_k\Big\{\mu_V\left(q_{ki}\rh_k\log(q_{ki}\rh_k)\right)\\
&\quad -\mu_V\left(q_{ki}\rh_k\right)\log\mu_V\left(q_{ki}\rh_k\right)+\mu_V\left(q_{ki}\rh_k\right)\log\mu_V(\rh_i)\Big\},~i\in\SS.
\end{align*} 
\end{lem}
\beg{proof}
Let $q_{ki}^{(n)}=q_{ki}\1_{[|x|\leq n]}$, $Q^{(n)}=(q^{(n)}_{ki})_{1\leq k,i\leq N}$ and $L^{[n]}=L^Z+Q^{(n)}$. Then   Lemma \ref{lem_mu} and Lemma \ref{wHar} imply that the semigroup generated by $L^{[n]}$ has an invariant measure with positive density $\rh^{[n]}$ w.r.t $\mu_{\pi,V}$ such that $\sup_{n\geq 1}\mu_{\pi,V}(\rh^{[n]}\log\rh^{[n]})<\infty$, which implies that $\rh_i^{[n]}\ra \rh_i$ weakly in $L^1(\mu_V)$.  Moreover, $\rh$ is a probability density and  also satisfies $\mu_V(\rh L g)=0$, $g\in C_c^2(\R^d\times\SS)$. Hence, Lemma \ref{wHar} implies that $\rh$ is positive.

The equation $\mu_{\pi,V}(\rh^{[n]}L^{[n]}g)=0$, $g\in C_c(\R^d\times\SS)$ yields, for all $i\in\SS$,
\beg{align}\label{equa_muLf}
&\pi_i\mu_V\left(\<Z_i,\nn f\>\rh_i^{[n]}\right)+\mu_V\left(\sum_{k=1}^N\pi_k q_{ki}^{[n]} f\rh_k^{[n]}\right)\nonumber\\
&=\pi_i\mu_V\left(\<\si^*_i\nn\rh_i^{[n]},\si^*_i\nn f\>\right), f\in C_c(\R^d).
\end{align}
Then there is $C_n>0$ which is independent of $f$ such that 
\beg{align}\label{rhH12}
\left|\mu_V\left(\<\si^*_i\nn\rh_i^{[n]},\si^*_i\nn f\>\right)\right|\leq C_n\sq{\mu_V(|\si_i^*\nn f|^2)+\mu_V(f^2)}
\end{align}
The defective log-Sobolev inequality yields the Poincar\'e inequality:
$$\mu_V(f^2)\leq C\mu_V(|\si_i^*\nn f|^2)+\mu_V(f)^2,~f\in H_{\si_i}^{1,2}(\R^d).$$
So the norm of $f$, $\sq{\mu_V(|\si_i^*\nn f|^2)+\mu_V(f)^2}$, which is induced by the following inner product of $g_1$, $g_2$ on $H_{\si_i}^{1,2}(\R^d)$:
$$\mu_V(g_1)\mu_V(g_2)+\mu_V(\<\si^*_i\nn g_1,\si^*_i\nn g_2\>),~g_1,~g_2\in H_{\si_i}^{1,2}(\R^d)$$
is  equivalent to the Sobolev norm on $H_{\si_i}^{1,2}(\R^d)$. Hence \eqref{rhH12} implies that $\rh_i^{[n]}\in H_{\si_i}^{1,2}(\R^d)$. Moreover, \eqref{equa_muLf} holds for all $f\in H_{\si_i}^{1,2}(\R^d)$. 

Next, we shall prove that $\sup_{n\geq 1}\mu_V(|\si_i^*\nn\sq{\rh_i^{[n]}}|^2)<\infty$. Let $f=\log(\rh_i^{[n]}+\ep)$ with $\ep\in (0,1)$. Then \eqref{equa_muLf} and $\sum_{j=1}^N\pi_j\mu_V(\rh_j^{[n]})=1$ imply that 
\beg{align*}
&2\pi_i\mu_V\left(|\si_i^*\nn\sq{\rh_i^{[n]}+\ep}|^2\right)+\pi_i\mu_V(q_{i}^{(n)}\rh_i^{[n]}\log(\rh_i^{(n)}+\ep))\\
&\leq \sum_{j\neq i}\pi_j\mu_V(q_{ji}^{(n)}\rh_j^{[n]}\log(\rh_i^{[n]}+\ep))+\ff {\pi_i} 2\mu_V(|Z_i|^2\rh_i^{[n]})\\
&\leq \sum_{j\neq i}\pi_j\Big\{\mu_V\left(q_{ji}^{(n)}\rh_j^{[n]}\log\left(q_{ji}^{(n)}\rh_j^{[n]}\right)\right)-\mu_V\left(q_{ji}^{(n)}\rh_j^{[n]}\right)\log\mu_V\left(q_{ji}^{(n)}\rh_j^{[n]}\right)\\
&\quad+\mu_V\left(q_{ji}^{(n)}\rh_j^{[n]}\right)\log\mu_V\left(\rh_i^{[n]}+\ep\right)\Big\}+\ff 1 2\left(\sup_{x\in\R^d}|Z_i|^2\right) \pi_i  \mu_V\left(\rh_i^{[n]}\right)\\
&\leq \sum_{j\neq i}\pi_j\Big\{\mu_V\left(q_{ji}^{(n)}\rh_j^{[n]}\log q_{ji}^{(n)}\right)+\mu_V\left(q_{ji}^{(n)}\rh_j^{[n]}\log \rh_j^{[n]}\right)+e^{-1}\\
&\quad+\ff {C_Q} {\pi_j}\log\left(\ff 1 {\pi_i}+1\right)\Big\}+\ff 1 2 \sup_{x\in\R^d}|Z_i|^2\\
&\leq \sum_{j\neq i}\Big\{\left(C_Q\log^+C_Q\right)+C_Q\pi_j\mu_V\left(\rh_j^{[n]}\log \rh_j^{[n]}\right)+\pi_j(C_Q+1)e^{-1}\\
&\quad+C_Q\log\left(\ff {\pi_i+1} {\pi_i}\right)\Big\}+\ff 1 2 \sup_{x\in\R^d}|Z_i|^2,
\end{align*}
and 
$$\pi_i\mu_V(q_{i}^{(n)}\rh_i^{[n]}\log(\rh_i^{[n]}+\ep))\geq \pi_i\mu_V(q_{i}^{(n)}\rh_i^{[n]}\log\rh_i^{[n]})\geq -\pi_i C_Q e^{-1}.$$
Hence, the monotone convergence theorem implies that there is $C>0$ which is independent of $i$ and $n$ such that 
$$\mu_V\left(\left|\si_i^*\nn\sq{\rh_i^{[n]}}\right|^2\right)=\lim_{\ep\ra 0^+}\mu_V\left(\left|\si_i^*\nn\sq{\rh_i^{[n]}+\ep}\right|^2\right)\leq C.$$
Since the defective log-Sobolev inequality implies the existence of a super Poincar\'e inequality, the essential spectrum of $L_i^0$ is empty, which implies that $H_{\si_i}^{1,2}(\R^d)$ is compactly embedded into $L^2(\mu_V)$. Thus,  $\sq{\rh_i}\in  H_{\si_i}^{1,2}(\R^d)$, and by a subsequence, $\rh_i^{[n]}$ converges to $\rh_i$ in $L^1(\mu_V)$ and $\sq{\rh_i^{[n]}}$ convergences to $\sq{\rh_i}$ weakly in $ H_{\si_i}^{1,2}(\R^d)$.

For $m>0$, let $f=\log\left((\rh_i^{[n]}\vee m)\we\ff 1 m\right)$. Then \eqref{equa_muLf} yields that
\beg{align}\label{inequ_sq}
&\pi_i\mu_V\left(\left|\si_i^*\nn\rh_i^{[n]}\right|^2\ff {\1_{[\ff 1 m\leq \rh_i^{[n]}\leq m]}} {\rh_i^{[n]}}\right)\nonumber\\
&=4\pi_i\mu_V\left(\left|\si_i^*\nn\sq{(\rh_i^{[n]}\vee m)\we \ff 1 m}\right|^2\right)\nonumber\\
&= \pi_i\mu_V\left(\<Z_i,\si_i^*\nn\rh_i^{[n]}\>\1_{[\ff 1 m\leq \rh_i^{[n]}\leq m]}\right)\nonumber\\
&\quad+ \sum_{k=1}^N\pi_k\mu_V\left( q_{ki}^{(n)}\left[\log\left((\rh_i^{[n]}\vee m)\we\ff 1 m\right)\right]\rh_k^{[n]}\right)\\
&\leq \ff {\pi_i} 2 \mu_V(|Z_i|^2\rh_i^{[n]})+\ff {\pi_i} 2 \mu_V\left(\left|\si_i^*\nn\rh_i^{[n]}\right|^2\ff {\1_{[\ff 1 m\leq \rh_i^{[n]}\leq m]}} {\rh_i^{[n]}}\right)\nonumber\\
&\quad +\sum_{k\neq i}\pi_k\mu_V\left(q_{ki}^{(n)}\rh_k^{[n]}\log\left((\rh_i^{[n]}\vee m)\we\ff 1 m\right)\right)\nonumber\\
&\quad -\pi_i\mu_V\left(q_i^{(n)}\rh_i^{[n]}\log\left((\rh_i^{[n]}\vee m)\we\ff 1 m\right)\right).\nonumber
\end{align}
Since $\rh_i^{[n]}\ra \rh_i$ in $L^1(\mu_V)$, $\sq{\rh_i^{[n]}}$ convergences to $\sq{\rh_i}$ weakly in $ H_{\si_i}^{1,2}(\R^d)$ and 
$$\sup_{n\geq 1,x\in\R^d}\left|q_i^{(n)}\log\left((\rh_i^{[n]}\vee m)\we\ff 1 m\right)\right|<\infty,$$
we have $\displaystyle{\lim_{n\ra\infty}\mu_V(|Z_i|^2\rh_i^{[n]})=\mu_V(|Z_i|^2\rh_i)}$,  
\beg{align*}
&\lim_{n\ra\infty}\sum_{k=1}^N\pi_k\mu_V\left(q_{ki}^{(n)}\rh_k^{[n]}\log\left((\rh_i^{[n]}\vee m)\we\ff 1 m\right)\right)\\
&=\sum_{k=1}^N\pi_k\mu_V\left(q_{ki}\rh_k\log\left((\rh_i\vee m)\we\ff 1 m\right)\right),
\end{align*}
and 
$$\varliminf_{n\ra\infty}\mu_V\left(\left|\si_i^*\nn\sq{(\rh_i^{[n]}\vee m)\we \ff 1 m}\right|^2\right)\geq \mu_V\left(\left|\si_i^*\nn\sq{(\rh_i\vee m)\we \ff 1 m}\right|^2\right).$$
Hence, combining these with \eqref{inequ_sq}, we obtain that
\beg{align}\label{inequ_muQ}
&2\pi_i\mu_V\left(\left|\si_i^*\nn\sq{(\rh_i\vee m)\we \ff 1 m}\right|^2\right)+\pi_i\mu_V\left(q_i\rh_i\log\left[(\rh_i\vee m)\we \ff 1 m\right]\right)\nonumber\\
&\leq \ff {\pi_i} 2 \mu_V(|Z_i|^2\rh_i)+\sum_{k\neq i}\mu_V\left(q_{ki}\rh_k\log\left((\rh_i\vee m)\we \ff 1 m\right)\right)\nonumber\\
&\leq \ff {\pi_i} 2 \mu_V(|Z_i|^2\rh_i)+\sum_{k\neq i}\Big\{\mu_V\left(q_{ki}\rh_k\log\left(q_{ki}\rh_k\right)\right)\\
&\quad -\mu_V\left(q_{ki}\rh_k\right)\log\mu_V\left(q_{ki}\rh_k\right)+\mu_V\left(q_{ki}\rh_k\right)\log\mu_V\left((\rh_i\vee m)\we \ff 1 m\right)\Big\}.\nonumber
\end{align}
Since $\rh_i>0$ and 
\beg{align*}
&\mu_V\left(\left|\si_i^*\nn\sq{(\rh_i\vee m)\we \ff 1 m}\right|^2\right)\\
&=\mu_V\left(\left|\si_i^*\nn\sq{\rh_i\vee m}\right|^2\right)-\mu_V\left(\left|\si_i^*\nn\sq{\rh_i}\right|^2\1_{[\rh_i\leq \ff 1 m]}\right),
\end{align*}
by monotone convergence theorem, we have
\beg{align*}
\varliminf_{m\ra\infty}\mu_V\left(\left|\si_i^*\nn\sq{(\rh_i\vee m)\we \ff 1 m}\right|^2\right)\geq \mu_V\left(\left|\si_i^*\nn\sq{\rh_i}\right|^2\right).
\end{align*}
Combining this with \eqref{inequ_muQ} and Fatou's Lemma, we have 
\beg{align*}
&\mu_V\left(\left|\si_i^*\nn\sq{\rh_i}\right|^2\right)+\pi_i\mu_V\left(q_i\rh_i\log\rh_i\right)\\
&\leq 2\pi_i\varliminf_{m\ra\infty}\mu_V\left(\left|\si_i^*\nn\sq{(\rh_i\vee m)\we \ff 1 m}\right|^2\right)\\
&\quad+\pi_i\varliminf_{m\ra\infty}\mu_V\left(q_i\rh_i\log\left[(\rh_i\vee m)\we \ff 1 m\right]\right)\\
&\leq \varliminf_{m\ra\infty}\Big\{2\pi_i\mu_V\left(\left|\si_i^*\nn\sq{(\rh_i\vee m)\we \ff 1 m}\right|^2\right)\\
&\quad+\pi_i\mu_V\left(q_i\rh_i\log\left[(\rh_i\vee m)\we \ff 1 m\right]\right)\Big\}\\
&\leq \ff {\pi_i} 2 \mu_V(|Z_i|^2\rh_i)+\sum_{k\neq i}\pi_k\Big\{\mu_V\left(q_{ki}\rh_k\log\left(q_{ki}\rh_k\right)\right)\\
&\quad -\mu_V\left(q_{ki}\rh_k\right)\log\mu_V\left(q_{ki}\rh_k\right)+\mu_V\left(q_{ki}\rh_k\right)\log\mu_V(\rh_i)\Big\}
\end{align*}

\end{proof}




The last result of this section is the following lemma on a priori estimate of entropy  of $\rh_i$ for $Z_i$ satisfies the integrability condition \eqref{ewZ}, which is crucial to the proof of Theorem \ref{thm_eu}.

\beg{lem}\label{entropy}

Let $\pi_{\min}=\displaystyle{\min_{1\leq k\leq N}\pi_k}$. Then, under the conditions of Lemma  \ref{ineq_nn_en} and Theorem \ref{thm_eu}, we have the following an a priori estimate
\beg{align}\label{inequ_ent}
\sum_{k=1}^N\pi_k\mu_V(\rh_k\log\rh_k)\leq &\max_{1\leq k\leq N}\left[\ff {v_k} {w_k}\mu_V(e^{w_k|\si_k^{-1}Z_k|^2})\right]-\log\pi_{\min}\nonumber\\
&+\tld C_Q \sum_{k=1}^N v_k+2\max_{1\leq k\leq N}\ff {v_k\be_k} {\ga_k},
\end{align}
where $\tld C_Q=C_Q\log^+ C_Q+2(C_Q+1)e^{-1}-2C_Q\log\pi_{\min}$. Moreover, there is a constant $C$ depending on $\pi_k,C_Q,w_k,v_k$ and $\mu_V(e^{w_k|Z_k|^2})$, $k\in \SS$ such that
\bequ\label{inequ_nn}
\sum_{k=1}^N\pi_k\mu_{V}\left(|\si_k^*\nn\sq {\rh_k}|\right)\leq C.
\enqu
\end{lem}
\beg{proof}
Let $\hat q_{ij}=\inf_{x\in\R^d} q_{ij}(x)$. Then by Lemma \ref{ineq_nn_en}, 
\beg{align}\label{nnleqent}
&2\pi_k\mu_V\left(|\si_k^*\nn\sq {\rh_k}|^2\right)+\pi_k\mu_V( q_k\rh_k\log\rh_k)\nonumber\\
&\leq \ff {\pi_k} {2w_k} \left(\mu_V\left(\rh_k\log\rh_k\right)-\mu_V(\rh_k)\log\mu_V(\rh_k)\right)\nonumber\\
&\quad +\ff {\pi_k} {2w_k}\mu_V(\rh_k)\log\mu_V\left(e^{w_k|Z_k|^2}\right)\\
&\quad+\sum_{j\neq k}\pi_j\Big(\bar q_{jk}\mu_V(\rh_j\log\rh_j)+\mu_V(\rh_j)\sup_{x\in\R^d}\left(q_{jk}(x)\log q_{jk}(x) \right)\nonumber\\
&\quad+(\bar q_{jk}-\hat q_{jk})e^{-1}+ e^{-1}-C_Q\mu_V(\rh_j)\log\pi_j\Big).\nonumber
\end{align}
Combining with (H\ref{H3}), we have
\beg{align}\label{entropy_1}
&\left(\ff {2\pi_k} {\ga_k}+\pi_k\bar q_k\right)\left(\mu_V(\rh_k\log\rh_k)-\mu_V(\rh_k)\log\mu_V(\rh_k)\right)-\ff {2\pi_k\be_k} {\ga_k}\mu_V(\rh_k)\nonumber\\
&\leq 2\pi_k\mu_V\left(|\si_k^*\nn\sq {\rh_k}|^2\right)+\pi_k\bar q_k\left(\mu_V(\rh_k\log\rh_k)-\mu_V(\rh_k)\log\mu_V(\rh_k)\right)\nonumber\\
&\leq \ff {\pi_k} {2w_k} \left(\mu_V\left(\rh_k\log\rh_k\right)-\mu_V(\rh_k)\log\mu_V(\rh_k)\right)\\
&\quad+\ff {\pi_k} {2w_k}\mu_V(\rh_k)\log\mu_V\left(e^{w_k|Z_k|^2}\right)\nonumber\\
&\quad+\sum_{j\neq k}\pi_j \bar q_{jk}\left(\mu_V(\rh_j\log\rh_j)-\mu_V(\rh_j)\log\mu_V(\rh_j)\right)\nonumber\\
&\quad+C_Q\log^+ C_Q+2(C_Q+1)e^{-1}-2C_Q\log\pi_{\min}.\nonumber
\end{align}

Let 
\beg{align*}
\ps&=(\ps_j)_{1\leq j\leq N}=\left(\mu_V(\rh_j\log\rh_j)-\mu_V(\rh_j)\log\mu_V(\rh_j)\right)_{1\leq j\leq N},\\
\Ps&=\ff {1} {2w_k}\mu_V(\rh_k)\log\mu_V\left(e^{w_k|Z_k|^2}\right).
\end{align*}
Then  \eqref{inequ_MM} and \eqref{entropy_1} yields that  
\beg{align*}
0&\leq \<\ps,(K+\bar Q)v\>_{\pi}+\<\Ps,v\>_\pi+\tld C_Q\sum_{k=1}^N v_k\leq -\<\ps,\1\>_{\pi}+\<\Ps,v\>_\pi+\tld C_Q\sum_{k=1}^N v_k,
\end{align*}
where $\<v,u\>_\pi=\sum_{i=1}^N\pi_iv_i u_i$, $u,v\in\R^N$. Hence
\beg{align*}
&\sum_{k=1}^N\pi_k\mu_V(\rh_k\log(\rh_k))\\
&\leq \<\Ps,v\>_\pi+\tld C_Q\sum_{k=1}^N v_k+\sum_{k=1}^N\pi_k\mu_V(\rh_k)\log\mu_V(\rh_k)+2\max_{1\leq k\leq N}\ff {\be_k v_k} {\ga_k}\\
&\leq  \<\Ps,v\>_\pi+ \tld C_Q\sum_{k=1}^N v_k-\log\pi_{\min}+2\max_{1\leq k\leq N}\ff {\be_k v_k} {\ga_k},
\end{align*}
in the last inequality, we use $\sum_{k=1}^N\pi_k\mu_V(\rh_k)=1$.
Then it is easy to see that the inequality \eqref{inequ_ent} holds. Consequently, \eqref{inequ_nn} holds by \eqref{nnleqent}. 

\end{proof}

\section{Proofs of main results}
\subsection{Th proof of Theorem \ref{thm_eu}}

Let $Z_i^{(n)}=Z_i\1_{[|Z_i|\leq n]}$, $\displaystyle\lim_{n\ra\infty}Z_i^{(n)}=Z_i$.  Then, according to Lemma \ref{lem_mu}, the Markov semigroup $P_t^{(n)}$ generated by $L^{(n)}$ has an invariant probability measure, denoted by $\rh^{(n)}\mu_{\pi,V}$. Moreover, Lemma \ref{entropy} yields that 
$$\sup_{n}\left(\mu_{\pi,V}(\rh^{(n)}\log\rh^{(n)})+\sum_{i=1}^N\mu_V\left(\pi_i|\si_i^*\nn\sq{\rh_i^{(n)}}|^2\right)\right)<\infty,$$
which implies there exists a $\rh\in L^1(\R^d\times\SS,\mu_{\pi,V})$ such that, by a subsequence, $\rh_i^{(n)}\ra \rh_i$ in $L^1(\mu_V)$ and $\sq{\rh_i^{(n)}}\ra \sq{\rh_i}$ weakly in $H_{\si_i}^{1,2}(\R^d)$.  Next, we shall prove that $\rh\mu_{\pi,V}$ is an invariant probability measure of $P_t$. Indeed, noticing that  
\beg{align*}
\mu_{\pi,V}(\rh^{(n)} f)&=\mu_{\pi,V}(\rh^{(n)} P_t^n f)\\
&=\mu_{\pi,V}(\rh^{(n)}(P_t^n f - P_t f))+\mu_{\pi,V}((\rh^{(n)}-\rh)P_tf)+\mu_{\pi,V}(\rh P_tf),
\end{align*}
and for $f\in\sB_b(\R^d\times\SS)$
\beg{align*}
\lim_{n\ra \infty}\mu_{\pi,V}((\rh^{(n)}-\rh)P_tf)=0,\qquad \lim_{n\ra\infty}\mu_{\pi,V}(\rh^{(n)} f)=\mu_{\pi,V}(\rh f),
\end{align*}
in order to prove  $\mu_{\pi,V}(\rh f)=\mu_{\pi,V}(\rh P_tf)$, we only need to prove that
$$\lim_{n\ra \infty}\mu_{\pi,V}(\rh_n(P_t^n f - P_t f))=0.$$
Let $(X_t,\La_t)$ be the process generated by $L^Q$. Then, by Lemma \ref{GirTh},  
\beg{align*}
&\mu_{\pi,V}(|P_t^n f - P_t f|)\\
&\leq \mu_{\pi,V}\Big(\E\Big|f(X_t,\La_t)\Big[e^{\int_0^t\left\<Z^{(n)}_{\La_s}(X_s),\d W_s\right\>-\ff 1 2\int_0^t|Z^{(n)}_{\La_s}(X_s)|^2\d s}\\
&\quad -e^{\int_0^t\left\<Z_{\La_s}(X_s),\d W_s\right\>-\ff 1 2\int_0^t|Z_{\La_s}(X_s)|^2\d s}\Big]\Big|\Big)\\
&\leq ||f||_\infty\left(\E e^{2\int_0^t\left\<Z_{\La_s^{\mu_{\pi,V}}}(X_s^{\mu_{\pi,V}}),\d W_s\right\>-\int_0^t\left|Z_{\La_s^{\mu_{\pi,V}}}(X_s^{\mu_{\pi,V}})\right|\d s}\right)^{\ff 1 2}\\
&\quad \times \left(\E\left|e^{\int_0^t\left\<(Z^{n}-Z)_{\La_s^{\mu_{\pi,V}}}(X_s^{\mu_{\pi,V}}),\d W_s\right\>-\ff 1 2\int_0^t\left|(Z^{n}-Z)\right|_{\La_s^{\mu_{\pi,V}}}^2(X_s^{\mu_{\pi,V}})\d s }-1\right|^2\right)^{\ff 1 2}.
\end{align*}
For all $\et>0$ and  all measurable function $h$ from $\R^d\times\SS$ to $\R^d$ with $|h(x,i)|\leq |Z_i(x)|$, $(x,i)\in \R^d\times\SS$, we have
\beg{align*}
&\E e^{2\et\int_0^t\<h(X_s^{\mu_{\pi,V}},\La_s^{\mu_{\pi,V}}),\d W_s\>-\et\int_0^t|h(X_s^{\mu_{\pi,V}},\La_s^{\mu_{\pi,V}})|}\\
&\leq \left(\E e^{4\et\int_0^t\<h(X_s^{\mu_{\pi,V}},\La_s^{\mu_{\pi,V}}),\d W_s\>-8\et^2\int_0^t|h(X_s^{\mu_{\pi,V}},\La_s^{\mu_{\pi,V}})|}\right)^{\ff 1 2}\\
&\quad \times\left(\E e^{(8\et^2-2\et)^+\int_0^t|h(X_s^{\mu_{\pi,V}},\La_s^{\mu_{\pi,V}})|^2\d s}\right)^{\ff 1 2}\\
&\leq \left(\ff 1 t\int_0^t\E e^{(8\et^2-2\et)^+t|h(X_s^{\mu_{\pi,V}},\La_s^{\mu_{\pi,V}})|^2}\d s\right)^{\ff 1 2}\\
&\leq \left(\ff 1 t\int_0^t\mu_{\pi,V}\left(P_s^Q e^{(8\et^2-2\et)^+t|h|^2}\right)\d s\right)^{\ff 1 2}\\
&\leq \left(\ff 1 t\int_0^te^{b_1s}\d s\right)^{\ff 1 2}\sq{\mu_{\pi,V}\left(e^{(8\et^2-2\et)^+t|Z|^2}\right)}\\
&< \infty,\qquad t<\ff {\min_{1\leq k\leq N} w_k} {(8\et^2-2\et)^+}.
\end{align*}
So, for $t$ small enough,  
\bequ\label{limPP}
\lim_{n\ra \infty}\mu_{\pi,V}(|P_t^n f - P_t f|)=0.
\enqu
Next, since $\{\rh^{(n)}\}_{n\geq 1}$ is uniformly integrable, for $\ep>0$, there is $m>0$ such that 
$\displaystyle{\sup_{n\geq 1}\mu_{\pi,V}(\rh^{(n)}\1_{[\rh^{(n)}\geq m]})\leq \ep}$. Thus
\beg{align*}
\mu_{\pi,V}(\rh^{(n)}|P_t^n f - P_t f|)&\leq \mu_{\pi,V}(\rh^{(n)}\1_{[\rh^{(n)}<m]}|P_t^n f - P_t f|)+2||f||_\infty\mu_{\pi,V}(\rh^{(n)}\1_{[\rh^{(n)}\geq m]})\\
&\leq m\mu_{\pi,V}(|P_t^n f - P_t f|)+2\ep ||f||_\infty,
\end{align*}
which combining with \eqref{limPP} yields that
$$\lim_{n\ra\infty}\mu_{\pi,V}(\rh^{(n)}|P_t^n f - P_t f|)=0.$$

\bigskip

Lastly, we shall prove our claim on the uniqueness. Since $||\si_i||$ is local bounded, it is easy to see from \eqref{ZA} that $\mu_V(\rh_i|\si_i Z_i|\1_{\cO})<\infty$ for all bounded domain $\cO\subset \R^d$. Then Lemma \ref{wHar} and Remark \ref{rem-add} yield that $\rh_i>0$, $i\in \SS$. By  \cite[Proposition 3.2.7, Proposition 3.2.5]{DPZ96}, there is only one invariant probability measure $\mu$ such that $\mu\ll\mu_{\pi,V}$.

\subsection{The proof of Theorem \ref{un_semigroup}}

Let 
\beg{align*}
\hat L_i^0=\ff 1 {\rh_i e^V}\div\left(e^V \rh_i a_i\nn \cdot\right),\qquad
\hat Z_i =Z_i-\si_i^*\nn \log\rh_i.
\end{align*}
Then $L_i^0=\hat L_i^0-a_i(\nn\log\rh_i)\cdot\nn$ and $ L^Z_i=\hat L_i^0+\si_i\hat Z_i\cdot\nn$. We shall prove first that  for $g_i\in L^\infty$, $i\in\SS$ so that the following equality holds for all  $f_i\in C_c^2(\R^d), i\in\SS$ 
\beg{align}\label{equ_mu_fg}
0=\sum_{i=1}^N\pi_i\mu_V\left(\left[g_i(1-\hat L_i^0)f_i  -\<\si_i\hat Z_i,\nn f_i\>g_i-\sum_{i=1}^N\pi_iq_{ij}f_j g_i\right]\rh_i\right),
\end{align}
then $g_i\in H^{1,2}_{loc}(\rh_i\d \mu_V)$. Indeed, fixing some $i\in\SS$, and letting  $f_k=0$, if $k\neq i$ and $f_i\in C_c^2(\R^d)$, we have
\beg{align}\label{equ_mu_fg-1}
\pi_i\mu_V\left(g_i \rh_i(1-\hat L_i^0)f_i \right)=\pi_i\mu_V\left(\<\si_i\hat Z_i,\nn f_i\>g_i\rh_i\right)+\sum_{k=1}^N\pi_k\mu_V\left(q_{ki}f_i g_k\rh_k\right).
\end{align} 
Let $\ze\in C_c^\infty(\R^d)$. Then
\beg{align*}
&\pi_i\mu_V\left((\ze g_i)\rh_i(1-\hat L_i^0)f_i \right)\\
&=\pi_i\mu_V\left(g_i\rh_i(1-\hat L_i^0)(\ze f_i) \right)+2\pi_i\mu_V\left(\<a_i\nn\ze,\nn f_i\>\rh_i\right)+\pi_i\mu_V\left((\hat L_i^0\ze)g_i\rh_i\right)\\
&=\pi_i\mu_V\left(\<\si_i\hat Z_i,\nn (\ze f_i)\>g_i\rh_i\right)+\sum_{k=1}^N\pi_k\mu_V\left(q_{ki}\ze f_i g_k\rh_k\right)\\
&\quad+2\pi_i\mu_V\left(\<a_i\nn\ze,\nn f_i\>g_i\rh_i\right)+\pi_i\mu_V\left((\hat L_i^0\ze)f_i g_i\rh_i\right).
\end{align*}
Since Lemma \ref{ex_den},  $\mu_V(e^{w_i|Z_i|^2})<\infty$ and that $||\si_i||$ is local bounded, there is a constant $C$ which is  independent  on $f_i$   such that
\beg{align*}
\left|\mu_V\left(\<\si_i\hat Z_i,\nn (\ze f_i)\>g_i\rh_i\right)\right|&\leq ||g_i||_{\infty}\sq{\mu_V\left((\ze+|\nn\ze|)^2|\hat Z_i|^2\rh_i\right)\mu_V\left((f_i^2+|\si_i^*\nn f_i|^2)\rh_i\right)}\\
&\leq C\sq{\mu_V\left((f_i^2+|\si_i^*\nn f_i|^2)\rh_i\right)}.
\end{align*}
By Lemma \ref{wHar}, we also  have 
\beg{align*}
\left|\sum_{k=1}^N\pi_k\mu_V\left(q_{ki}\ze f_i g_k\rh_k\right)\right|
&\leq \sq{\mu_V(f_i^2\rh_i)\left|\mu_V\left(\left(\sum_{k=1}^N\pi_kq_{ki} g_k\rh_k\right)^2\ff {\ze^2} {\rh_i}\right)\right|}\\
&\leq \ff {\sq{\mu_V(f_i^2\rh_i)}} {\inf_{x\in\supp \ze} \rh_i(x)}\sq{\left|\mu_V\left(\left(\sum_{k=1}^N\pi_kq_{ki} g_k\rh_k\right)^2\ze^2\right)\right|}\\
&\leq C\sq{\mu_V(f_i^2\rh_i)},
\end{align*}
for some positive constant  $C$ which is independent of $f_i$. Similarly, 
$$2\mu_V\left(\<a_i\nn\ze,\nn f_i\>\rh_i\right)+\mu_V\left((\hat L_i^0\ze)g_i f_i\rh_i\right)\leq C\sq{\mu_V(f_i^2\rh_i)}.$$
Hence
$$\left|\mu_V\left((\ze g_i)\rh_i(1-\hat L_i^0)f_i \right)\right|\leq C\sq{\mu_V((f_i^2+|\si_i^*\nn f_i|^2)\rh_i)},~f_i\in C_c^2(\R^d),$$
which implies that 
$$\ze g_i\in H^{1,2}_{\si,0}(\rh_i\d\mu_V),$$
where $H^{1,2}_{\si_i,0}(\rh_i\d\mu_V)$ is the complete of $C_c^2(\R^d)$ under the following norm of $f$
$$\sq{\mu_V((f^2+|\si_i^*\nn f|^2)\rh_i)}.$$
Thus $g_i\in H^{1,2}_{loc}(\rh_i\d \mu_V)$. 

Next, we shall prove that there is only one extension of $(L, C_c^2(\R^d\times\SS))$ that generates a  $C_0$-semigroup in $L^1(\R^d\times\SS,\mu)$. To this end, we only have to prove that $(1-L)(C_c^2(\R^d\times\SS))$ is dense in $L^1(\R^d\times\SS,\mu)$. Let $\ze_n\in C_c^\infty(\R^d)$ with $\1_{B_n(0)}\leq \ze_n\leq \1_{B_{2n}(0)}$ and $|||\nn\ze_n|||_\infty\leq \ff 1 n$. Let $g_i\in L^\infty$, $i\in\SS$ such that \eqref{equ_mu_fg} hold. Fixing $i\in\SS$, and letting $f_k=0$ if $k\neq i$; $f_i=\ze^2_n g_i$, we have, following from \eqref{equ_mu_fg},  
\beg{align}\label{equa_zeg}
0&=\pi_i\mu_V(\ze^2_n g^2_i\rh_i)+\pi_i\mu_V\left(\<a_i\nn(\ze_n^2 g_i),\nn g_i\>\rh_i\right)\nonumber\\
&\quad -\pi_i\mu_V\left(\<\si_i\hat Z_i,\nn(\ze_n^2 g_i)\>g_i\rh_i\right)-\mu_V\left(\sum_{k=1}^N\pi_k q_{ki}(\ze_n^2 g_i)g_k\rh_k\right)\nonumber\\
&=\pi_i\mu_V(\ze_n^2 g^2_i)+\pi_i\mu_V\left(\<a_i\nn(\ze_n g_i),\nn (\ze_n g_i)\>\rh_i\right)-\pi_i\mu_V\left(\<a_i\nn\ze_n,\nn \ze_n\> g_i^2\rh_i\right)\\
&\quad -\pi_i\mu_V\left(\<\si_i\hat Z_i,\nn(\ze_n^2 g_i)\>g_i\rh_i\right)-\mu_V\left(\sum_{k=1}^N \pi_k  q_{ki}(\ze_n^2 g_i)g_k\rh_k\right).\nonumber
\end{align}
Letting 
$$f_k(x)=(\ze_n g_i)^2(x)\1_{[k=i]}(k),~x\in\R^d,~k\in\SS,$$ since $(\ze_n g_i)^2\in H_{\si_i}^{1,2}(\rh_i\d \mu_V)$, we have $\mu(\hat L^0f)=0$ and $\mu(Lf)=0$. Then the equality $\mu(\hat L^0f)=\mu(Lf)$ yields that  
\beg{align*}
0&=\pi_i\mu_V\left(\<\si\hat Z_i,\nn(\ze_n^2 g_i^2)\>\rh_i\right)+\mu_V\left(\sum_{k=1}^N\pi_kq_{ki}\ze_n^2 g_i^2\rh_k\right)\\
&=\pi_i\mu_V\left(\<\si\hat Z_i,\nn(\ze_n^2 g_i)\>g_i\rh_i\right)+\pi_i\mu_V\left(\<\hat Z_i, \nn g_i\>g_i\ze_n^2\rh_i\right)+\mu_V\left(\sum_{k=1}^N\pi_kq_{ki}\ze_n^2 g_i^2\rh_k\right)\\
&=2\pi_i\mu_V\left(\<\si\hat Z_i, \nn g_i\>g_i\ze_n^2\rh_i\right)+\pi_i\mu_V\left(\<\si\hat Z_i,\nn \ze_n^2\>g_i^2\rh_i\right)+\mu_V\left(\sum_{k=1}^N\pi_kq_{ki}\ze_n^2 g_i^2\rh_k\right).
\end{align*}
Thus 
\beg{align*}
\pi_i\mu_V\left(\<\si\hat Z_i, \nn g_i\>g_i\ze_n^2\rh_i\right)
&=-\pi_i\mu_V\left(\<\hat Z_i,\nn \ze_n\>\ze_n g_i^2\rh_i\right)-\ff 1 2\mu_V\left(\sum_{k=1}^N\pi_kq_{ki}\ze_n^2 g_i^2\rh_k\right),\\
\pi_i\mu_V\left(\<\si_i\hat Z_i,\nn(\ze_n^2 g_i)\>g_i\rh_i\right)&=2\pi_i\mu_V\left(\<\si_i\hat Z_i,\nn \ze_n\>\ze_n g_i^2\rh_i\right)+\pi_i\mu_V\left(\<\si_i\hat Z_i,\nn g_i\>g_i\ze_n^2\rh_i\right)\\
&=\pi_i\mu_V\left(\<\si_i\hat Z_i,\nn \ze_n\>\ze_n g_i^2\rh_i\right)-\ff 1 2\mu_V\left(\sum_{k=1}^N\pi_kq_{ki}(\ze_n g_i)^2\rh_k\right).
\end{align*}
Combining with \eqref{equa_zeg}, we arrive at 
\beg{align}\label{zeqq}
&\pi_i\mu_V(\ze_n^2 g^2_i\rh_i)+\pi_i\mu_V\left(\<a_i\nn(\ze_n g_i),\nn (\ze_n g_i)\>\rh_i\right)\nonumber\\
&=\pi_i\mu_V\left(\<a_i\nn\ze_n,\nn \ze\> g_i^2\rh_i\right)+\mu_V\left(\sum_{k=1}^Nq_{ki}(\ze_n^2 g_i)g_k\rh_k\right)\\
&\quad+\pi_i\mu_V\left(\<\si_i\hat Z_i,\nn \ze_n\>\ze_n g_i^2\rh_i\right)-\ff 1 2\mu_V\left(\sum_{k=1}^N\pi_kq_{ki}(\ze_n g_i)^2\rh_k\right).\nonumber
\end{align}
Since
\beg{align*}
&\sum_{k=1}^N\pi_kq_{ki}g_ig_k\rh_k-\ff 1 2 \sum_{k=1}^N\pi_kq_{ki} g_i^2\rh_k\\
&=-\ff {\pi_i} 2  q_i g_i^2\rh_i+\ff 1 2 \sum_{k\neq i}^N\pi_kq_{ki}g_i( g_k- g_i)\rh_k+\ff 1 2\sum_{k\neq i}^N\pi_k q_{ki}g_ig_k\rh_k\\
&=-\ff {\pi_i} 2  q_i g_i^2\rh_i-\ff 1 2 \sum_{k\neq i}^N\pi_k q_{ki}( g_k- g_i)^2\rh_k\\
&\quad+\ff 1 2 \sum_{k\neq i}^N\pi_k q_{ki}g_k( g_k- g_i)\rh_k+\ff 1 2\sum_{k\neq i}^N\pi_k q_{ki}g_ig_k\rh_k\\
&=-\ff {\pi_i} 2  q_i g_i^2\rh_i-\ff 1 2 \sum_{k\neq i}^N\pi_k q_{ki}( g_k- g_i)^2\rh_k+\ff 1 2 \sum_{k\neq i}^N\pi_kq_{ki}g_k^2\rh_k,
\end{align*}
we have that
\beg{align}\label{qqgg}
&\sum_{i=1}^N\left(\sum_{k=1}^N\pi_kq_{ki}g_ig_k\rh_k-\ff 1 2 \sum_{k=1}^N\pi_kq_{ki} g_i^2\rh_k\right)\nonumber\\
&=-\ff 1 2 \sum_{i=1}^N\pi_i  q_i g_i^2\rh_i -\ff 1 2\sum_{i=1}^N\sum_{k\neq i}^N\pi_kq_{ki}( g_k- g_i)^2\rh_k+\ff 1 2 \sum_{i=1}^N\sum_{k\neq i}^N\pi_k q_{ki}g_k^2\rh_k\nonumber\\
&=-\ff 1 2\sum_{i=1}^N\sum_{k\neq i}^N\pi_kq_{ki}( g_k- g_i)^2\rh_k -\ff 1 2 \sum_{i=1}^N\pi_i  q_i g_i^2\rh_i+\ff 1 2 \sum_{k=1}^N\sum_{i\neq k}^N\pi_kq_{ki}g_k^2\rh_k\\
&=-\ff 1 2\sum_{i=1}^N\sum_{k\neq i}^N\pi_k q_{ki}( g_k- g_i)^2\rh_k. \nonumber
\end{align}
Hence, \eqref{zeqq}, \eqref{qqgg} and the definition of $\ze_n$ yield that
\beg{align*}
&\sum_{i=1}^N\left\{\pi_i\mu_V\left(\ze_n^2 g^2_i \rh_i+\<a_i\nn(\ze_n g_i),\nn (\ze_n g_i)\>\rh_i\right)+\mu_V\left(\ff {\ze_n^2} 2\sum_{k\neq i}^N\pi_k q_{ki}( g_k- g_i)^2\rh_k\right)\right\}\\
&=\sum_{i=1}^N\pi_i\mu_V\left(\<a_i\nn\ze_n,\nn \ze_n\> g_i^2\rh_i\right)+\sum_{i=1}^N\pi_i\mu_V\left(\<\si_i\hat Z_i,\nn \ze_n\>\ze_n g_i^2\rh_i\right)\\
&\leq \ff {\max_{i\in\SS}||g_i||^2_\infty} {n^2} \sum_{i=1}^N\mu_V(||a_i||\rh_i)+\ff {\max_{i\in\SS}||g_i||^2_\infty} n\sum_{i=1}^N\pi_i\mu_V\left(|\si_i\hat Z_i|\rh_i\right).
\end{align*}
For $r$ small enough, we have 
\beg{align}\label{ZA}
&\mu_V(||a_i||\rh_i)+\mu_V\left(|\si_i\hat Z_i|\rh_i\right)\nonumber\\
&\leq \ff 1 2\mu_V\left(\rh_i(|Z_i|^2+3||\si_i||^2)\right)+\mu_V(||\si_i||\cdot|\si_i^*\nn \rh_i|)\nonumber\\
&\leq \ff 1 2 \mu_V\left(\rh_i|Z_i|^2\right)+2\mu_V\left(\rh_i||\si_i||^2\right)+2\mu_V\left(|\si_i^*\nn\sq \rh_i|\right)\nonumber\\
&\leq  \ff 1 {2r} \mu_V(\rh_i\log\rh_i)+2\mu_V\left(|\si_i^*\nn\sq \rh_i|\right)-\ff 1 {2r}\mu_V(\rh_i)\log\mu_V(\rh_i)\\
&\quad+\ff 1 {2r}\Big\{\log\mu_V\left(e^{r|Z_i|^2}\right)+\log\mu_V\left(e^{3r||\si_i||^2}\right)\Big\}\nonumber\\
&<\infty,~\qquad i\in\SS.\nonumber
\end{align}
Therefore
$$\varlimsup_{n\ra\infty}\sum_{i=1}^N\pi_i\mu_V(\ze_n^2 g^2_i)\leq 0,$$
which implies that $g_i=0$, $i\in\SS$.




\end{document}